\documentclass[preprint,12pt]{elsarticle}

\usepackage{amssymb}
\usepackage{amsmath,amssymb}
\usepackage{tabularx}
\usepackage{hyperref}
 \usepackage{booktabs}
\newcolumntype{Y}{>{\hsize=1.5\hsize}X}
\newcolumntype{B}{X}
\newcolumntype{T}{>{\hsize=0.5\hsize}X}
\newcolumntype{Z}{>{\hsize=1.8\hsize}X}
\usepackage{xcolor}
\usepackage{lineno}
\usepackage[mathscr]{euscript}
\usepackage{stmaryrd}
\usepackage{mathtools}
\usepackage[symbol]{footmisc}
\journal{Physics Reports}

\begin{document}

\begin{frontmatter}



\title{Stability of Ecological Systems: A Theoretical Review}


\author{Can Chen$^{a,b}$, Xu-Wen Wang$^b$, and Yang-Yu Liu$^{b,c,}$\footnote[1]{Corresponding author. Email: yyl@channing.harvard.edu}}

\affiliation{organization={School of Data Science and Society and Department of Mathematics, University of North Carolina at Chapel Hill},
            city={Chapel Hill},
            postcode={27599}, 
            state={NC},
            country={USA}}

\affiliation{organization={Channing Division of Network Medicine, Department of Medicine, Brigham and Women's Hospital, Harvard Medical School},
            city={Boston},
            postcode={02115}, 
            state={MA},
            country={USA}}

\affiliation{organization={Carl R. Woese Institute for Genomic Biology, Center for Artificial Intelligence and Modeling, University of Illinois at Urbana-Champaign},
            city={Champaign},
            postcode={61801}, 
            state={IL},
            country={USA}}

\begin{abstract}
The stability of  ecological systems is a fundamental concept in ecology, which offers profound insights into species coexistence, biodiversity, and community persistence. In this article, we provide a systematic and comprehensive review on the theoretical frameworks for analyzing the stability of  ecological systems. Notably, we survey various stability notions,  including linear stability, sign stability, diagonal stability, D-stability, total stability, sector stability, structural stability, and higher-order stability. For each of these stability notions, we  examine   necessary or sufficient conditions for achieving such stability and demonstrate the intricate interplay of these conditions on the network structures of ecological systems. Finally, we explore the future prospects of these stability notions.

\end{abstract}




\end{frontmatter}
\tableofcontents

\section{Introduction}\label{sec:1}
Ecological systems are usually subject to continual disturbances or perturbations, and their responses are often characterized quantitatively as stability, a fundamental concept rooted in dynamical systems and control theory \cite{buma2020disturbances,kefi2019advancing,li2006global,neubert1997alternatives,pascual2005criticality,saunders1975stability,schoon2012understanding,vsiljak1974connective,townley2007predicting}. 
Loosely speaking, an ecological system is stable if its state represented by a vector of its species abundances, do not change too much under small perturbations.  
Therefore,
understanding the stability of ecological systems is particularly important due to its implications for species coexistence, biodiversity, and community persistence \cite{cottingham2001biodiversity,goh1975stability, grimm1992application, harrison2001dynamical,holling1973resilience,levine2017beyond,ozaki2020understanding,tokeshi2009species}.  More importantly, these implications can contribute to the development of effective interventions and policies that maintain the overall  health of ecological systems \cite{domptail2013managing,francis2021management,goh2012management,janssen2006governing}.

Numerous mathematicians, physicists, and ecologists have explored the stability of  ecological systems. In his pioneer work, May linearized  an unspecified nonlinear dynamical system around an equilibrium  to perform the local linear stability analysis of ecological systems \cite{May-Nature-1972, May-Book-1974,may1975stability}.
The local linear stability  is completely characterized by the eigenvalues of the Jacobian matrix (of the unspecified dynamical model), which is often called the community matrix in the ecological contexts. If all the eigenvalues of the community matrix have negative real parts,  the ecological system is locally stable around the equilibrium. By incorporating random matrix theory, May successfully derived a simple condition to characterize the stability of  ecological systems. While there have been intense debates about the relationship between stability and complexity over the past five decades \cite{mccann2000diversity, odenbaugh2011complex,pimm1984complexity,roberts1974stability,vsiljak1975complex,upadhyay2000stability,veloz2020complexity,yonatan2022complexity}, May's result offers a powerful foundation for analyzing the stability of ecological systems with diverse structural characteristics, such as interaction types \cite{Allesina-Nature-2012}, degree heterogeneity \cite{yan2017degree}, self-regulation \cite{barabas2017self}, modularity \cite{grilli2016modularity}, and more.

May's framework does not need to know the dynamical model of the ecological system. Instead, it focuses on the Jacobian matrix of the unspecified dynamical model. Hence, the framework is model-implicit. In the past, actually many ecological models, such as the generalized Lotka-Volterra (GLV) model (1910) \cite{bhargava1989generalized,malcai2002theoretical},  Holling type II model (1959) \cite{agiza2009chaotic,skalski2001functional}, and  MacArthur's consumer-resource model (1970) \cite{advani2018statistical,chesson1990macarthur}, have been proposed to study the dynamics of ecological systems. Remarkably, the GLV model is the most commonly used model in the stability analysis of  ecological systems due to its simplicity \cite{al2003stability,altieri2021properties,beretta1988generalization,beretta1988global,bunin2017ecological,goh1977global, kon2011age,zeeman1998three}. Most of these findings rely on Lyapunov theory, which allows us to determine the stability of a system without explicitly integrating the differential equation. Lyapunov theory encompasses two fundamental methods. One is based on linearization (e.g., May's approach), which is often referred to as the Lyapunov indirect method, while the other is called the Lyapunov direct method. The essence of Lyapunov direct method is to construct the so-called Lyapunov function, an “energy-like” scalar function whose time variation can be viewed as the ``energy dissipation.'' It is not restricted to small perturbations, and in principle can be applied to all dynamical systems. However, we lack a general theory to find a suitable Lyapunov function for an arbitrary dynamical system. Instead, we have to rely on our experience and intuition. Notably, certain diagonal-type Lyapunov functions are  useful for the stability analysis of  ecological systems described by the GLV model \cite{Kaszkurewicz-Book-00}.

\begin{figure}[t]
    \centering
    \includegraphics[scale=0.34]{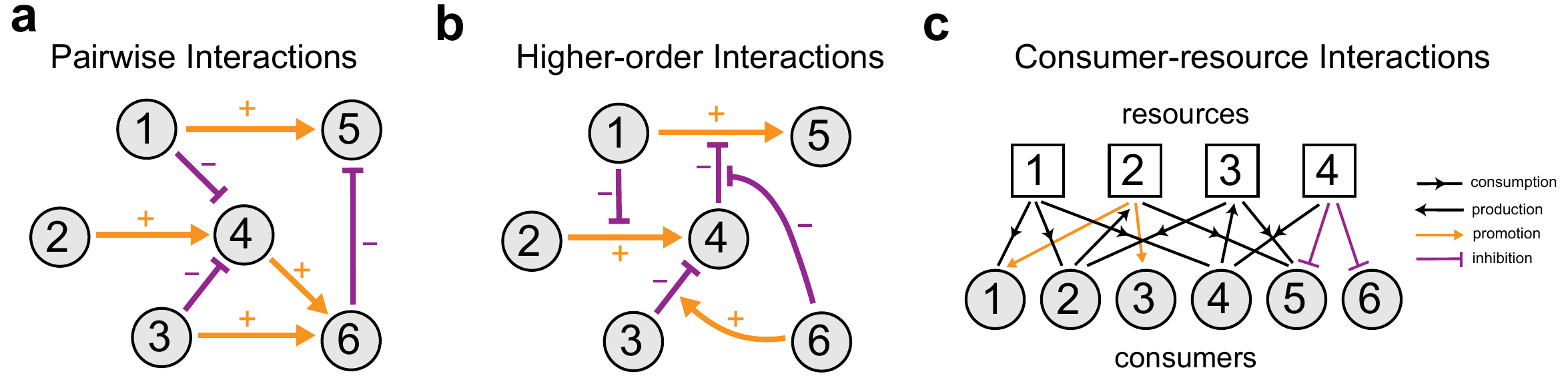}
    \caption{Pairwise interactions versus higher-order interactions in  ecological systems. (a) Pairwise interspecies  interactions. (b) Higher-order interspecies interactions. (c) Consumer-resource interactions naturally implies higher-order interspecies  interactions.}
    \label{fig:hypergraph}
\end{figure}

The feasibility of an ecological system, which concerns whether a specific set of conditions or parameters can support the system over time, can also impact its stability \cite{aladwani2020ecological,aparicio2023feasibility,goh1977feasibility,grilli2017feasibility,paulus1995feasibility,stone2018feasibility}. Intriguingly, the relationship between feasibility and stability is noteworthy when dealing with ecological systems described by the GLV model or  MacArthur's consumer-resource model. It has been proved that under certain conditions, the feasibility of such systems can imply their  stability \cite{case1979global,stone2018feasibility}. Moreover, feasibility is  closely related to a concept known as structural stability, which reflects the ability of an ecological system to qualitatively maintain its dynamics and overall behavior under small disturbances or perturbations of the system model itself \cite{rohr2014structural}. A great amount of metrics have been proposed to quantify the structural stability of  ecological systems through feasibility analysis  \cite{grilli2017feasibility,rohr2014structural,saavedra2014structurally}.

Recently, there has been a growing interest in higher-order interactions within  ecological systems, which refers to the intricate relationships that extend beyond pairwise interspecies interactions  \cite{bairey2016high,barbosa2023experimental,gibbs2022coexistence,grilli2017higher, li2020advances}, see Figure \ref{fig:hypergraph}a, b. Mathematically, these higher-order interactions emerged in an ecological system can be naturally represented as a hypergraph, where its hyperedges can connect an arbitrary number of nodes \cite{berge1984hypergraphs}. The resulting dynamics can be expressed in the form of polynomials \cite{chen2021controllability}. The analysis of hypergraph dynamics often involves the application of tensor theory \cite{chen2020tensor,chen2021controllability,pickard2023hat, surana2022hypergraph}, which deals with multidimensional arrays generalized from vectors and matrices \cite{chen2022explicit,chen2019multilinear,chen2021multilinear}. Lyapunov theory can also be leveraged to analyze the stability of polynomial dynamical systems \cite{ahmadi2019algebraic}. Furthermore, the consumer-resource model can be implicitly considered as an instance of higher-order interactions, with consumer species as nodes and resources as hyperedges, see Figure \ref{fig:hypergraph}c. Increasing evidence has revealed that higher-order interactions play a significant role in the stability of  ecological systems \cite{bairey2016high,gibbs2022coexistence,grilli2017higher}.

\begin{table}[t!]
\centering
\caption{Summary of stability notions that have been used to study  ecological systems.}
\begin{tabularx}{\textwidth}{|T|X|Y|X|}
\hline
Notion  & Main Model  & Characteristics &  Key References \\
\hline
\scriptsize{Linear \newline Stability} &  \scriptsize{Linear model (\ref{eq: linearized}), \newline GLV model (\ref{eq: glv})}  & \scriptsize{Characterizing the eigenvalue spectrum of the community matrix \textbf{M} through random matrix theory.}  & \scriptsize{\cite{Allesina-Nature-2012,allesina2015predicting,barabas2017self,baron2020dispersal,baron2023breakdown,galla2018dynamically,gibbs2018effect,gravel2016stability, grilli2016modularity,jacquet2016no,krumbeck2021fluctuation,May-Nature-1972, May-Book-1974, pigani2022delay, poley2023generalized, sidhom2020ecological, stone2018feasibility, tang2014correlation, Yan-arXiv-14, yan2017degree, yang2023time}
}
\\\hline

\scriptsize{Sign \newline Stability} & \scriptsize{Linear model (\ref{eq: linearized})} & \scriptsize{Assessing stability  based solely on the signs of the elements in the community matrix \textbf{M} without considering their magnitudes.}  &  \scriptsize{\cite{allesina2008network,dambacher2003qualitative,  haraldsson2020model, Jeffries-74, logofet1982sign, Maybee-69, quirk1965qualitative, Yamada-Networks-90}
}\\\hline

\scriptsize{Diagonal Stability} & \scriptsize{GLV model (\ref{eq: glv})} & 
\scriptsize{The diagonal stability of the interaction matrix \textbf{A}, i.e., there exists a diagonal matrix \textbf{P} such that $\textbf{A}^\top\textbf{P}+\textbf{P}\textbf{A}\prec 0$,  implies the Lyapunov stability of the GLV model.}
& \scriptsize{\cite{arcak2011diagonal,arcak2008passivity, arcak2006diagonal, 4434115, beretta1988generalization, beretta1988global, gibbs2018effect, goh1977global, grilli2017feasibility, hernandez1997lotka, hou2015global, Kaszkurewicz-Book-00, kon2011age, motee2012stability, Redheffer-SIAM-85,simpson2021diagonal, wang2012diagonal,worz1978global}} \\\hline

\scriptsize{D-Stability} & \scriptsize{GLV model (\ref{eq: glv})} &
\scriptsize{The D-stability of the interaction matrix \textbf{A}, i.e., the matrix \textbf{X}\textbf{A} is stable for any positive diagonal matrix \textbf{X}, implies that any feasible equilibrium of the GLV model is locally stable.}  &  \scriptsize{\cite{arrow1958note, chu2007equivalent,enthoven1956theorem,gibbs2018effect,johnson1974sufficient,rohr2014structural}}\\\hline 

\scriptsize{Total \newline Stability} & \scriptsize{GLV model (\ref{eq: glv})} & 
\scriptsize{The total stability of the interaction matrix \textbf{A}, i.e., any principal sub-matrix of \textbf{A} is D-stable, implies that a new feasible equilibrium of the GLV model remains locally stable after the removal or extinction of certain species.}
 &  \scriptsize{\cite{logofet2018matrices,logofet2005stronger,logofet1987on,quirk1965qualitative}}\\\hline 

 \scriptsize{Sector \newline Stability} & \scriptsize{Generic population dynamics model (\ref{eq:gpm})} & 
\scriptsize{\textcolor{black}{A semi-feasible equilibrium point is called sector stable if every solution that starts within a non-negative neighborhood remains in the same or an even larger non-negative neighborhood and eventually converges to that equilibrium point.}}
 &  \scriptsize{\cite{case1979global,goh1978sector}}\\\hline 

\scriptsize{Structural Stability} & \scriptsize{GLV model (\ref{eq: glv}) } &\scriptsize{Measuring the capacity to qualitatively maintains the dynamics under small disturbances or perturbations through feasibility analysis of the GLV model.} &  \scriptsize{\cite{Andronov-1937,grilli2017feasibility,Peixoto-1962,pettersson2020predicting,recknagel1985analysis,rohr2014structural,saavedra2017structural,saavedra2014structurally,Smale-1961,Smale-1967,song2018structural}} \\\hline 

\scriptsize{Higher-order Stability} & \scriptsize{GLV model with \newline higher-order interactions  (\ref{eq:hoglv}), \newline MacArthur’s consumer-resource mode (\ref{eq:cr})} & \scriptsize{Analyzing the stability of ecological systems with (implicit) higher-order interactions via polynomial systems theory and local stability analysis.} & \scriptsize{\cite{ahmadi2019algebraic,akimenko2021stability,aladwani2019addition,aparicio2023feasibility,bairey2016high,bieg2023stability,butler2018stability,case1979global,chen2022explicit,de2000random,din2021discrete,duffy2016identifying,el2023age,gawronski2022instability,gibbs2022coexistence,grilli2017higher,kelsic2015counteraction,schreiber1996global,singh2021higher,synodinos2023effects, yoshino2008rank}} \\
\hline
\end{tabularx}
\label{tab:0}
\end{table}

Most significantly, the stability of an ecological system is profoundly determined and characterized by its underlying network structure, which captures how species influence one another and how disturbances or perturbations propagate among species. In this article, we provide a systematic and inclusive literature review on the theoretical frameworks for analyzing the stability of  ecological systems. 
In  Section \ref{sec:2}, we briefly review the fundamentals of Lyapunov  theory, which acts as a foundation for the stability analysis of ecological systems. From  Section \ref{sec:3} to  \ref{sec:10}, we comprehensively survey various stability notions, including linear stability, sign stability, diagonal stability, D-stability, total stability, sector stability, structural stability, and higher-order stability, see Table \ref{tab:0}. For each of these stability notions, we examine necessary or sufficient conditions that are required to establish such stability and illustrate the complicated relationships between these conditions and the network structures of ecological systems. Finally, we discuss the future prospects of these stability notions in  Section \ref{sec:11} and conclude in  Section \ref{sec:12}.

\section{Preliminaries}\label{sec:2}
\textcolor{black}{In this section, we  review the fundamentals of Lyapunov  theory, which plays a significant role in the various stability analyses of  ecological systems. We adapt most of the notations and results from the comprehensive work of \cite{khalil2009lyapunov,khalil2015nonlinear,Murray-Book-Robotics,sastry1999lyapunov,Slotine-Book-91}.}

\subsection{Definitions of Lyapunov stability}\label{sec:2.1}
Consider a general $S$-dimensional nonlinear system 
\begin{equation}
    \dot{\bf x}(t) = {\bf f}({\bf x}(t), t) \quad \text{with} \quad {\bf
  x}(t_0)={\bf x}_0, \quad {\bf x}(t) \in \mathbb{R}^S.
\label{eq:GNL}
\end{equation}
In most cases, we are interested in the stability of the system near an equilibrium point ${\bf x}^* \in \mathbb{R}^S$ that satisfies ${\bf f}({\bf x}^*, t)= {\bf 0}$. Without loss of generality, we can shift the origin of the system and assume that the equilibrium point of interest occurs at ${\bf
x}^*={\bf 0}$. Let $\textbf{M}$ be the Jacobian matrix evaluated at $\textbf{x}^*$. Denote $n_{-}$, $n_0$, and $n_{+}$ be the numbers  of eigenvalues (counting multiplicities) of $\textbf{M}$ with negative, zero, and positive real part, respectively. An equilibrium point is called hyperbolic if $n_0 = 0$, i.e., there are no eigenvalues on the imaginary axis. A hyperbolic equilibrium point is called a hyperbolic saddle if $n_{-} n_{+} \neq 0$. Since a generic matrix has $n_0 = 0$,  equilibrium points in a generic dynamical system are typically hyperbolic. 

\begin{figure}[t!]
\centering
\includegraphics[width=0.6\textwidth]{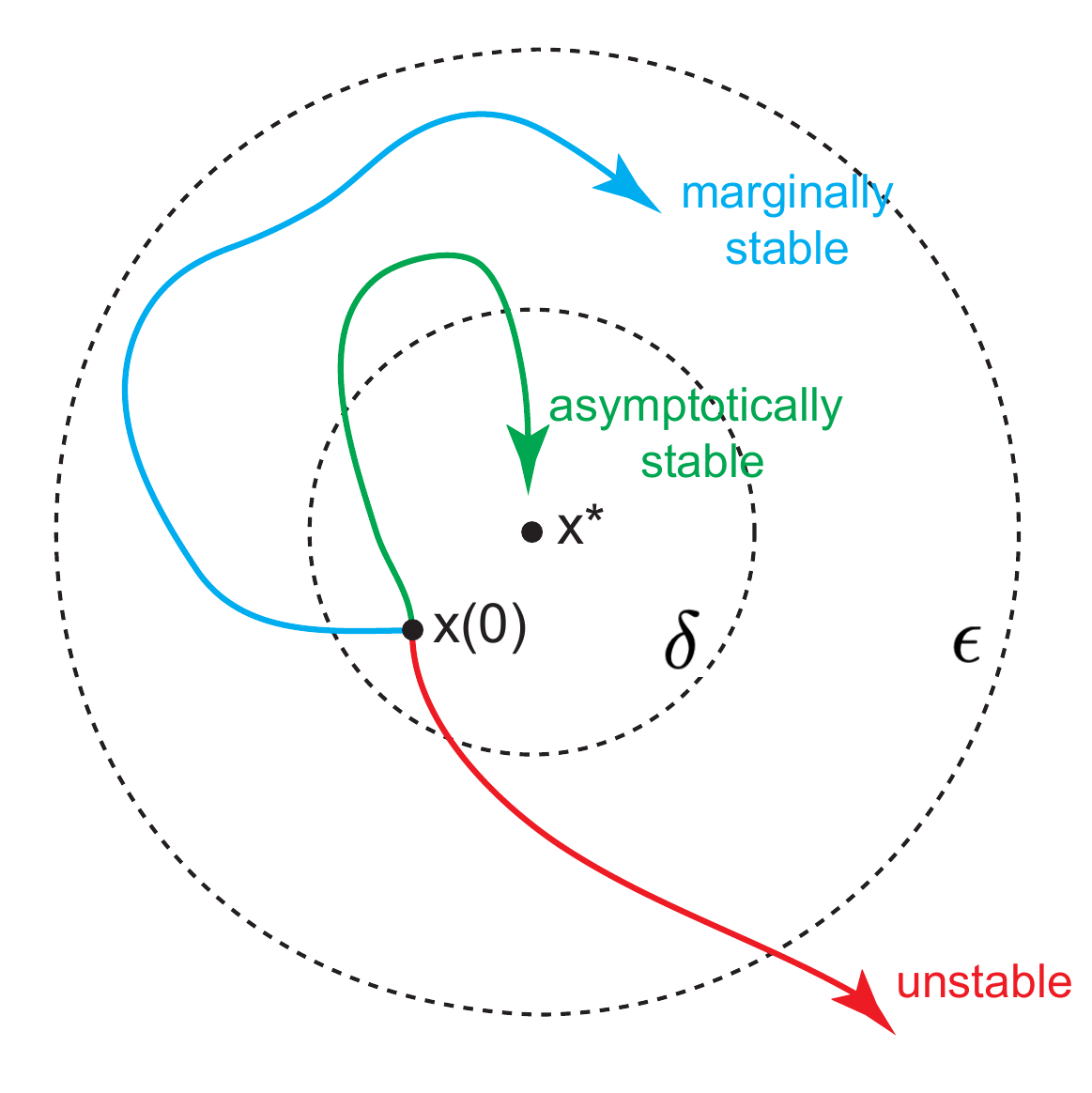}
\caption{Geometrical implication of stability. 
For nonlinear systems $\dot{\textbf{x}}(t) = \textbf{f} (\textbf{x}(t), t)$, due to the complex and rich behavior of nonlinear dynamics, various types of stability, e.g., stability, asymptotic stability, and global asymptotic stability, can be
discussed. Intuitively and roughly speaking, if all solutions of the system that start out near an equilibrium point ${\bf x}^*$ stay near ${\bf x}^*$ forever, then ${\bf x}^*$ is stable (in the sense of Lyapunov). More strongly, if ${\bf x}^*$ is stable and all solutions that start out near ${\bf x}^*$ converge to ${\bf x}^*$ as $t \to \infty$, then ${\bf x}^*$ is asymptotically stable. ${\bf x}^*$ is marginally stable if it is stable but not asymptotically stable. This figure was redrawn from \cite{Slotine-Book-91}. 
}
\label{fig:stability-demo}
\end{figure}

An equilibrium point ${\bf x}^*={\bf 0}$ is stable (in the sense of Lyapunov) at $t=t_0$ if $\forall \epsilon > 0$,  $\exists\, \delta(\epsilon, t_0)>0$ such that 
\begin{equation*}
    ||{\bf x}(t_0)|| \le \delta(\epsilon, t_0) \Rightarrow ||{\bf x}(t)|| \le \epsilon, \text{ } \forall t \ge t_0 > 0.
\end{equation*}
Otherwise, ${\bf x}^*$ is unstable, see Figure \ref{fig:stability-demo}. An equilibrium point ${\bf x}^*={\bf 0}$ is asymptotically stable at $t=t_0$ if it is stable and locally attractive, i.e., $\exists\, \delta(t_0)>0$ such that 
\begin{equation*}
    ||{\bf x}(t_0)|| \le \delta(t_0) \Rightarrow \lim_{t\to \infty} {\bf x}(t)= {\bf 0}.
\end{equation*}
An equilibrium point which is stable but not asymptotically stable is called marginally stable, see Figure~\ref{fig:stability-demo}. Both stability and asymptotic stability are defined at a time instant $t_0$. In practice, it is often desirable for a system to have a certain uniformity in its behavior. Uniform (asymptotic) stability requires that the equilibrium point is (asymptotically) stable for all $t_0>0$. Note that notions of uniformity are only relevant for time-varying or non-autonomous systems. For autonomous or time-invariant systems, (asymptotic) stability naturally implies uniform (asymptotic) stability. The definition of asymptotic stability does not quantify the rate of convergence. The notion of exponential stability guarantees a minimal rate of convergence. An equilibrium point ${\bf x}^*={\bf 0}$ is exponentially stable if $\exists\, \delta, \alpha, \lambda >0$ such that
\begin{equation*}
    ||{\bf x}(t_0)|| \le \delta \Rightarrow ||{\bf x}(t)|| \le \alpha ||{\bf x}(t_0)|| e^{-\lambda (t-t_0)}, \text{ }\forall t \ge t_0.
\end{equation*}
Exponential stability is a very strong form of stability because it implies uniform asymptotic stability.

The above definitions of stability, asymptotic stability, and exponential stability are local and only describe the behavior of a system near an equilibrium point. We call an equilibrium point globally (asymptotically or exponentially) stable, if it is (asymptotically or exponentially) stable for all initial conditions ${\bf x}_0 \in \mathbb{R}^S$. Though global stability is very desirable, it is very difficult to achieve in many systems.  Note that linear time-invariant (LTI) systems are either asymptotically stable, or marginally stable, or unstable. Moreover, linear asymptotic stability is always global and exponential, and linear instability
always implies exponential blow-up. The above refined notions of stability are explicitly needed only for nonlinear  systems. 

\subsection{Lyapunov’s indirect method}\label{sec:2.2}
Lyapunov's indirect  method is based on linearization, and is  concerned with the local stability of a nonlinear  system and is based on the intuition that a nonlinear  system should behave similarly to its linearized approximation for small disturbances or perturbations. The linearization can be used to give a conservative bound on the domain of attraction of the equilibrium point for the original nonlinear system. 

For a general nonlinear system (\ref{eq:GNL}) with continuously differentiable ${\bf f}({\bf x}(t), t)$ around the equilibrium point ${\bf x}^*={\bf 0}$,  the system dynamics can be rewritten as
\begin{equation*}
    \dot{\bf x}(t) = {\bf M}(t) {\bf x}(t) + {\bf f}_\text{h.o.t}({\bf x}(t), t),
\end{equation*}
where ${\bf M}(t) = \frac{\partial {\bf f}({\bf x}, t)}{\partial {\bf x}}|_{{\bf x}={\bf x}^*}$ is the Jacobian matrix of ${\bf f}({\bf x},t)$ with respect to ${\bf x}$, evaluated at the equilibrium point ${\bf x}^*$, and ${\bf f}_\text{h.o.t}({\bf x}, t)$ represents the higher-order terms in ${\bf x}$. We typically require ${\bf f}_\text{h.o.t}({\bf x}, t)$ approaches zero uniformly, which is obviously true for an autonomous system. Then $\dot{\bf x}(t) = {\bf M}(t) {\bf x}(t)$ is called the uniform linearization of the original nonlinear system (\ref{eq:GNL}) at the equilibrium point ${\bf x}^*={\bf 0}$. If this uniform linearization exists and the Jacobian matrix ${\bf M}(t)$ is bounded, then the uniform asymptotic stability of the equilibrium point ${\bf x}^*={\bf 0}$ for the linearization implies its uniform local asymptotic stability for the original nonlinear system. 

For an autonomous system $\dot{\bf x}(t)={\bf f}( \textbf{x}(t))$, the linearization is simply an LTI system $\dot{\bf x}(t)={\bf M} {\bf x}(t)$. Denote the 
eigenvalues of ${\bf M}$ as $\lambda({\bf M})$. We have the following results: (i) if the linearization is strictly stable (i.e., all the eigenvalues of $\textbf{M}$ lie in the closed left half of the complex plane), then the equilibrium point ${\bf x}^*={\bf 0}$ is asymptotically stable for the original nonlinear system  (note that a real square matrix $\textbf{M}$ that satisfies $\text{Re}[\lambda({\bf M})] < 0$  is often called Hurwitz stable); (ii) if the linearization is unstable (i.e., at least one eigenvalue of ${\bf M}$ has positive real part), then the equilibrium point is unstable for the original nonlinear system; (iii) if the linearization is marginally stable (i.e, all eigenvalues of ${\bf M}$ are in the left-half complex plane, but at least one of them is on the imaginary axis), then one can not conclude anything from the linearization---the equilibrium point may be stable, asymptotically stable, or unstable for the original nonlinear system. 

\subsection{Lyapunov’s direct method}\label{sec:2.3}
Lyapunov's direct method  (also called the second method of Lyapunov) can be considered as a mathematical extension of a fundamental physical observation. If the total energy of a mechanical or electrical system is continuously dissipated, then the system must eventually settle down to an equilibrium point, no matter whether it is linear or nonlinear. Thus, we may conclude the stability of a system by studying the energy change rate of the system,
without explicitly solving the original differential equation (\ref{eq:GNL}).

Let ${\bf B}_\epsilon$ be a ball of size $\epsilon$ around the origin ${\bf x}={\bf 0}$, i.e., ${\bf B}_\epsilon = \{ {\bf x} \in \mathbb{R}^S: ||{\bf x}|| < \epsilon\}$. A scalar continuous function $V({\bf x}) : \mathbb{R}^S \to \mathbb{R}$ is locally positive definite (LPD) if  $V({\bf 0})=0$ and ${\bf x} \neq {\bf 0} \Rightarrow V({\bf x}) > 0$ in the ball ${\bf B}_\epsilon$. If $V({\bf 0})=0$ and the above property holds for the whole state space, then $V({\bf x})$ is globally positive definite (GPD). Similarly, a scalar continuous function $V({\bf x})$ is locally (or globally) positive semi-definite if $V({\bf 0}) = 0$ and  ${\bf x} \neq {\bf 0} \Rightarrow V({\bf x}) \ge 0$ in the ball ${\bf B}_\epsilon$ (or for the whole state space). A scalar continuous function $V({\bf x}, t): \mathbb{R}^S \times \mathbb{R}_{>0} \to \mathbb{R}$ is locally positive definite if $V({\bf 0}, t)=0$ and there exists a time-invariant LPD function $V_0({\bf x})$ that is dominated by $V({\bf x}, t)$, i.e., 
$V({\bf x}, t) \ge V_0({\bf x}), \text{ } \forall t \ge t_0$. If $V({\bf 0},t )=0$ and the above property holds for the whole state space, then $V({\bf x}, t)$ is globally positive definite. A scalar continuous function $V({\bf x}, t)$ is called decrescent if $V({\bf 0}, t)=0$ and there exists a time-invariant LPD function $V_1({\bf x})$ that dominates $V({\bf x}, t)$, i.e., $V({\bf x}, t) \le V_1({\bf x}), \text{ }\forall t \ge t_0$. A scalar continuous function $V({\bf x}, t)$ is called radially unbounded if $\lim_{||x||\to \infty} V({\bf x}, t) \to \infty$ uniformly on $t$. 

Let $V({\bf x}, t)$ be a non-negative function with derivative $\dot{V}({\bf x}, t)$ along the trajectories of the system, i.e., 
\begin{equation*}
    \dot{V}({\bf x}, t)=\frac{\partial V}{\partial t} + \frac{\partial V}{\partial {\bf x}} {\bf f}({\bf x}, t).
\end{equation*}
Lyapunov's theorems of equilibrium point stability are summarized in Table  \ref{tab:Lyapunov}. Note that they are all sufficiency theorems. If for a particular choice of Lyapunov function candidate $V$, the condition on $\dot{V}$ is not met, we cannot draw any conclusions on the system's stability. Many Lyapunov functions may exist for the same system. Specific choices of Lyapunov functions may yield more precise results than others.

\begin{table}[t]
\caption{Summary of Lyapunov stability theorems on equilibrium point ${\bf x}^* ={\bf 0}$ \cite{Murray-Book-Robotics}.}
\label{tab:Lyapunov}
\begin{tabularx}{\textwidth}{|X|T|Y|}
  \hline                   
$V({\bf x}, t)$  &  $-\dot{V}({\bf x}, t)$
& Stability (in the sense of Lyapunov) \\
\hline
LPD                   & LPSD & stable \\ \hline
LPD, decrescent & LPSD & uniformly stable \\\hline
LPD, decrescent & LPD & uniformly asymptotically stable\\ \hline
GPD, decrescent, \newline radially unbounded & GPD & globally uniformly asymptotically \newline stable  \\
  \hline  
\end{tabularx}
\end{table}

For an LTI system $\dot{\bf x}(t) = {\bf M} {\bf x}(t)$, we often use a quadratic Lyapunov function candidate $ V= \textbf{x}^\top \textbf{P} \textbf{x}$, where $\textbf{P}$ is a symmetric positive definite matrix (i.e., $\textbf{x}^\top\textbf{P} \textbf{x} >0, \text{ }\forall \textbf{x} \neq {\bf 0}$), denoted as $\textbf{P}\succ 0$. It is easy to derive that
\begin{equation*}
    \dot{V}= \dot{\textbf{x}}^\top \textbf{P} \textbf{x} + \textbf{x}^\top \textbf{P} \dot{\textbf{x}} = -\textbf{x}^\top \textbf{Q} \textbf{x},
\end{equation*}
where 
\begin{equation}
{\bf Q} = {\bf Q}^\top = - (\textbf{M}^\top \textbf{P} + \textbf{P} \textbf{M}). \label{eq:Lyapunovequation}
\end{equation}
The above equation  is often called the Lyapunov equation of the LTI system. If $\textbf{Q}$ is positive definite, then the system is globally asymptotically stable. However, if $\textbf{Q}$ is not positive definite, then no stability conclusion can be drawn. Fortunately, Lyapunov proved a necessary and sufficient condition for an LTI system to be strictly stable. For any symmetric positive definite matrix $\textbf{Q}$, the unique matrix solution $\textbf{P}$ of the Lyapunov equation (\ref{eq:Lyapunovequation}) is symmetric positive definite. Therefore, we can start by choosing a simple positive definite matrix $\textbf{Q}$ (e.g., the identity matrix $\textbf{I}$), then solve for $\textbf{P}$ from the Lyapunov equation, and finally verify whether $\textbf{P}$ is positive definite.

\section{Linear stability analysis}\label{sec:3}
As a simple application of Lyapunov's indirect method, linear stability analysis has been extensively applied to  ecological  systems, helping us better understand the intricate relationships between stability and biodiversity \cite{Allesina-Nature-2012,allesina2015stability,barclay1978deterministic,hogg1989stability,landi2018complexity,May-Nature-1972,May-Book-1974,nunney1980stability,weinstock1987local}. \textcolor{black}{More significantly, under some special conditions, local stability determined by linear stability analysis implies global stability for certain ecological systems \cite{hsu1995global}.} 

In his pioneer work, May considered an ecological system of $S$ species coexisting at  a feasible equilibrium point ${\bf x}^*$ of a dynamical system $\dot{\bf x}(t) = {\bf f}({\bf x}(t))$ that describes the time-dependent abundance vector ${\bf x}(t)$ of the $S$ species \cite{May-Nature-1972,May-Book-1974}. 
Suppose that one species is subjected to a small but sudden population increase or decrease. As a result, other species populations may show immediate changes away from the equilibrium.
The manipulated species itself, if with initially increased (or reduced) abundance, may begin to decline (or increase) toward the equilibrium because of self-regulation. These immediate, direct changes can be referred to as first-order effects and  described by the linearized equation of the original dynamical system, i.e., 
\begin{equation}\label{eq: linearized}
    \dot{\textbf{z}}(t) = \textbf{M}\textbf{z}(t),
\end{equation}
where ${\bf z}(t)={\bf x}(t)-{\bf x}^*\in\mathbb{R}^S$ denotes the deviation from the equilibrium, and the Jacobian matrix  ${\bf M} = \frac{\partial{\bf f}({\bf x})}{\partial{\bf x}}|_{{\bf x}={\bf x}^*}\in\mathbb{R}^{S\times S} $ is often referred to as the  community matrix. The diagonal elements of the community matrix $\textbf{M}_{ii} = \frac{\partial f_i({\bf x})}{\partial x_i} |_{\textbf{x}={\bf x}^*}$ represent the self-regulation of species $i$, while  the off-diagonal elements $\textbf{M}_{ij} = \frac{\partial f_i({\bf x})}{\partial x_j} |_{\textbf{x}={\bf x}^*} $ capture the impact that species $j$ has on species $i$ around the equilibrium point ${\bf x}^*$.

Numerous studies have explored the stability of the linearized equation (\ref{eq: linearized}), most of which employ random matrix theory to characterize the eigenvalue spectrum of the community matrix \cite{Allesina-Nature-2012,allesina2015predicting, baron2020dispersal, baron2023breakdown, gibbs2018effect,May-Nature-1972}. These studies can be broadly classified into two categories. The first category, referred to as model-implicit approaches, focuses on analyzing the stability of the linearized equation without requiring specific knowledge about the underlying dynamics $\textbf{f}(\textbf{x})$. The second category, referred to as model-explicit approaches, incorporates additional information about $\textbf{f}(\textbf{x})$, e.g., the GLV model, to provide insights into the stability analysis.

\subsection{Model-implicit approaches}\label{sec:3.1}
Since the empirical parameterization of the exact functional form of ${\bf f}({\bf x})$ is  difficult for an  ecological system, model-implicit approaches directly utilize the linearized representation to perform stability analysis. 

\subsubsection{May's classical result}\label{sec:3.1.1}
May  considered that $\textbf{M}_{ij}$ are randomly drawn from a distribution with mean $\mu=0$ and variance $\sigma^2$ with probability $C$ and are 0 otherwise \cite{May-Nature-1972}. Hence, $\sigma$ represents the characteristic interspecies interaction strength and $C$ is the ratio between actual and potential interactions in the ecological system (often referred to as the connectance). For simplicity, the diagonal elements are chosen to be the same with $-d=-1$, representing the intrinsic damping timescale of each species, so that if disturbed from equilibrium, it would return with such a damping time by itself. May found that for random interactions drawn from a Gaussian distribution $\mathcal{N}(0, \sigma^2)$, a randomly assembled system is stable (in the sense that all the eigenvalues of the community matrix  have negative real parts) if the so-called ``complexity'' measure
\begin{equation}
    \sigma \sqrt{C S} < d=1.
\end{equation}
This implies that more complexity (i.e., larger $C S$) tend to destabilize community dynamics \cite{May-Nature-1972,May-Book-1974}.

\begin{table}[t]
 \caption{Stability criteria for  ecological systems with different interaction types. In each case, the stability criterion is derived for a large community matrix ${\bf M}$. Diagonal elements $\textbf{M}_{ii}$, representing self-regulation, are set to be $-d$. For random interactions, off-diagonal elements $\textbf{M}_{ij}$ are randomly drawn from a normal distribution $\mathcal{N}(0,\sigma^2)$ with probability $C$ and are 0 otherwise. Predator-prey interactions come in pairs with opposite signs, i.e., $\textbf{M}_{ij}>0$ and $\textbf{M}_{ji}<0$.  With probability $C$, we sample $\textbf{M}_{ij}$ from $|\mathcal{N}(0,\sigma^2)|$ and $\textbf{M}_{ji}$ from $-|\mathcal{N}(0,\sigma^2)|$. With probability $(1-C)$, both $\textbf{M}_{ij}$ and  $\textbf{M}_{ji}$ are set to be 0. Community matrices with mutualistic and/or competitive interactions can be constructed similarly \cite{Allesina-Nature-2012}.}
\label{tab:Allesina}
\centering
  \begin{tabularx}{\textwidth}{|X|X|}
\hline
    Interaction type & Stability Criterion \\
\hline
    Random & $\sigma \sqrt{SC}<d$ \\
\hline
Random with  correlation ($\rho$) & $\sigma \sqrt{SC}(1+\rho)<d$ \\
\hline
Random with degree heterogeneity ($\xi$) & $\sigma \sqrt{(S-1)C\xi}<d$ \\
\hline
    Predator-Prey & $\sigma \sqrt{SC}<d \frac{\pi}{\pi-2}$\\
\hline
    Mixture of mutualism & $\sigma \sqrt{SC}<d \frac{\pi}{\pi+2}$ \\
and competition & \\
\hline
    Mutualism & $\sigma (S-1)C \sqrt{\frac{2}{\pi}} <d$ \\
\hline
    Competition & $\sigma \left[ \sqrt{S C}
      \frac{\pi+2-4C}{\sqrt{\pi(\pi-2C)}}  + C \sqrt{\frac{2}{\pi}} \right] <d$ \\
\hline
  \end{tabularx}
\end{table}

May's result is closely related to Girko's Circular Law  \cite{girko2004strong} in random matrix theory. Consider a random $S \times S$ real matrix ${\bf M}$ with entries independent and taken randomly from a normal distribution $\mathcal{N}(0,\sigma^2)$. Then as $S\to \infty$, the eigenvalues of ${\bf M}/\sqrt{S \sigma^2}$ are uniformly distributed in the unit disk
centered at $(0,0)$ in the complex plane \cite{Girko-TPA-1984}. Sommers et al. considered possible correlations between the off-diagonal elements $\textbf{M}_{ij}$ and $\textbf{M}_{ji}$ and 
proved that the eigenvalues of \textbf{M} are uniformly distributed in an ellipse with the real and imaginary directions $1+\rho$ and $1-\rho$ (where $\rho=\mathbb{E}[\textbf{M}_{ij}\textbf{M}_{ji}]/\text{Var}(\textbf{M}_{ij})$ for $i\neq j$), respectively, known as Sommers' Elliptical Law \cite{Sommers-PRL-1988}. 

\begin{figure}[th!]
\includegraphics[scale=0.48]{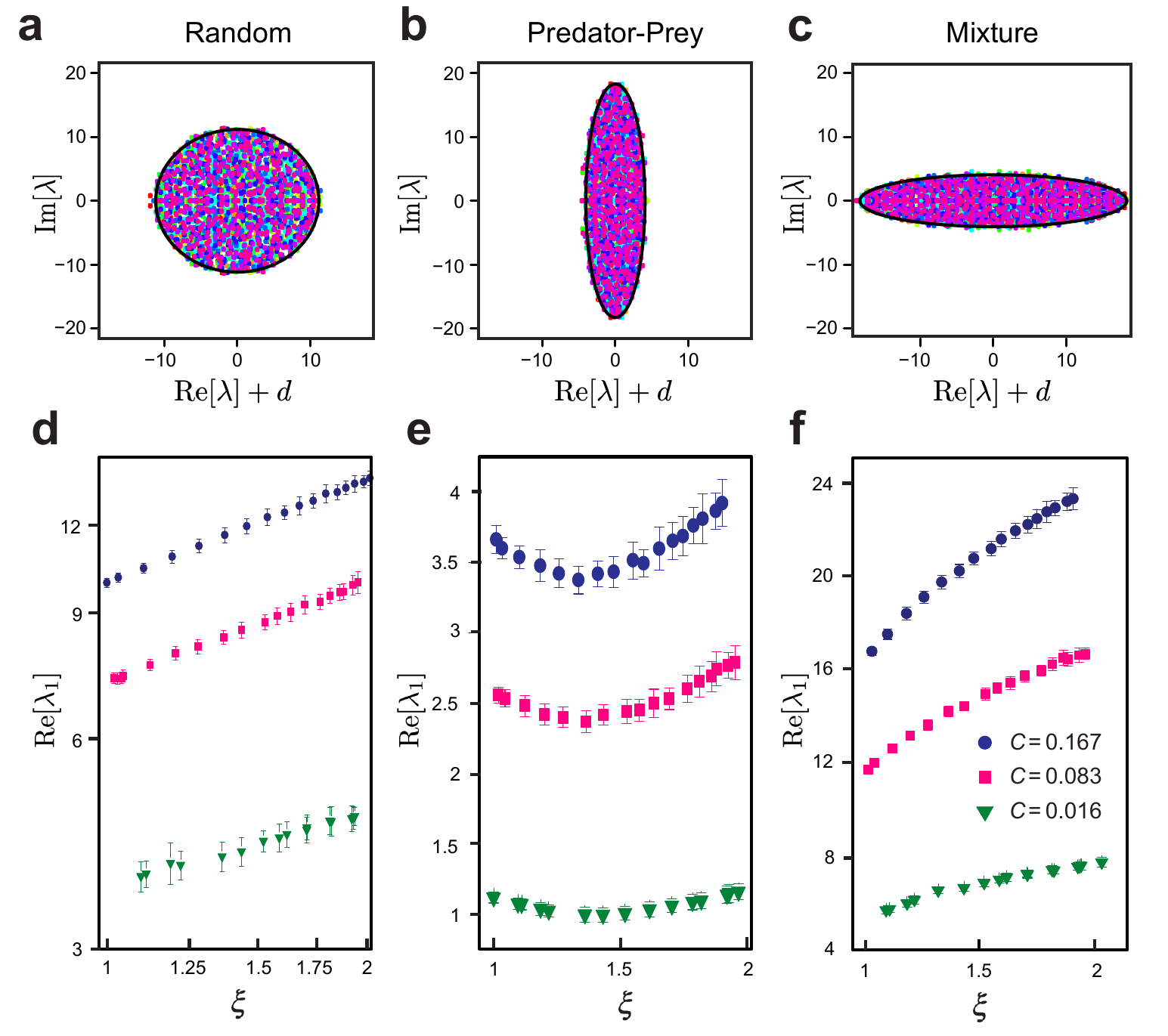}
  \caption{Eigenvalue distributions of 10 community matrices (color) with $-d=-1$ on the diagonal and off-diagonal elements following the random (a), predator-prey (b), and mixture of competition and mutualism interactions (c), respectively. $S=250$, $C=0.25$, $\sigma=1$. The black ellipses are analytical results.  Panels a, b, and c were redrawn from  \cite{Allesina-Nature-2012}. \textcolor{black}{Impact of degree heterogeneity on the stability of ecological systems with random interactions (d), predator-prey (e), and mixture of competition and mutualism interactions (f). The dots represent the results from numerical simulations on 3-modal networks with different connectances. Higher $\text{Re}[\lambda_1]$ indicates lower stability. Each error bar represents the standard deviation of 100 independent runs. $S=1200$, $\sigma=1$, and $d=0$. Panel d was drawn in the log-log scale. Panels d, e, and f were redrawn from \cite{yan2017degree}.}}
  \label{fig: eig}
\end{figure}

May's result continues to be influential almost five decades later, not because it asserts that  ecological systems must be inherently unstable, but rather because it highlights the importance of specific structural characteristics that enable real ecological systems to maintain stability despite their inherent complexity \cite{Bascompte-Science-2010}. In other words, nature must adopt some devious and delicate strategies to cope with this stability-complexity paradox.

\subsubsection{Impacts of interaction types and correlation}\label{sec:3.1.2}
One of the specificity is the existence of well-defined interspecific relationships observed in nature, e.g., predator-prey, competition, mutualism, and a mixture of  competition and mutualism. In 2012, by leveraging Sommers' Ellipical Law \cite{Sommers-PRL-1988}, Allesina et al. refined May's result and provided  stability criteria for all these interspecific interaction types  \cite{Allesina-Nature-2012} , as shown in  Table \ref{tab:Allesina} and Figure \ref{fig: eig}a, b, and c. They found remarkable differences between predator-prey interactions, which are stabilizing, but mutualistic and competitive interactions, which are destabilizing. Additionally, the correlation between interaction strengths determines interaction types. Tang et al. therefore incorporated  correlation  into the stability analysis by deriving a new stability criterion for large  ecological systems with random interactions \cite{tang2014correlation} , i.e., 
\begin{equation}\label{eq:corr}
    \sigma\sqrt{SC}(1+\rho)-\mu<d,
\end{equation}
where $\rho$ and $\mu$ are the overall correlation between pairs of interactions and the mean of the off-diagonal elements in the community matrix $\textbf{M}$, respectively. \textcolor{black}{The criterion (\ref{eq:corr}) can be viewed as a generalization of those for predator-prey, competition, and mutualism, derived by Allesina et al. \cite{Allesina-Nature-2012}.}
The authors found that the effect of correlation of interaction strengths substantially influences the stability of large food webs (with predator-prey interactions) compared to other network structural properties. Hence, the presence of correlation between interactions of species can significantly influence the locations of the eigenvalues of the community matrix and the resulting stability of ecological systems \cite{Allesina-Nature-2012,baron2022eigenvalues, kuczala2016eigenvalue,tang2014correlation}.

\subsubsection{Impacts of degree heterogeneity}\label{sec:3.1.3}
Many of the ``devious strategies'' adopted by nature can now be tested with the revised formula as a reference point. For example, one can further study the impact of degree heterogeneity on the stability of  ecological systems. Degree heterogeneity measures the variability of the number of interactions associated with each species. 
Yan et al.  found that for ecological systems with random interactions or a mixture of competition and mutualism interactions, increasing the degree of heterogeneity always destabilizes ecological systems \cite{yan2017degree}, see Figure \ref{fig: eig}d, f.
For ecological systems with predator-prey interactions  constructed from either simple network models or a realistic food web model (cascade model), high heterogeneity is always destabilizing, yet moderate heterogeneity is stabilizing \cite{Yan-arXiv-14,yan2017degree}, see Figure \ref{fig: eig}e. 
Furthermore, they obtained a stability criterion for large  ecological systems with random interactions by considering the factor of degree heterogeneity, i.e., 
\begin{equation}
    \sigma\sqrt{(S-1)C\xi}<d,
\end{equation}
where $\xi$ denotes the degree heterogeneity. Similarly, Allesina et al.  found that  broad degree distributions tend to stabilize food webs by approximating the real part of the leading (``rightmost'') eigenvalue of the community matrix \cite{allesina2015predicting}. Recently, Baron derived a closed-form expression for the eigenvalue spectrum of a general directed and weighted network \cite{baron2022eigenvalue}. The findings of this study suggest  that network heterogeneity appears to be a destabilizing influence in most circumstances, except when the interactions are very asymmetric (e.g., very asymmetric predator-prey interactions). These results are consistent with what were reported by Yan et al. \cite{yan2017degree}.

\subsubsection{Impacts of self-regulation}\label{sec:3.1.4}
Self-regulation (reflected as negative diagonal elements of the community matrix)  is also a key factor for stability in real or random ecological systems. To study the effect of self-regulation on the stability of ecological systems, Barabás et al.  derived an analytic approximation of these diagonal elements, i.e., 
\begin{equation}\label{eq:3.5}
    d_{ii}=-d=-2^{1.5}\text{Re}[\lambda_1]
\end{equation}
for $i=1,2,\dots, S$, where $\lambda_{1}$ is the  leading (``rightmost'') eigenvalue of the community matrix \textbf{M} without self-regulations \cite{barabas2017self}. Notably, empirical food webs can only achieve stability when the majority of species exhibit strong self-regulation \cite{barabas2017self}. Even for random ecological systems, attaining stability also requires negative self-regulations for a large amount of species \cite{barabas2017self}. In addition, Tang et al. investigated the effect of non-constant diagonal on the eigenvalue distribution of the community matrix \textbf{M} \cite{tang2014correlation}. They found that a moderate variance of the diagonal elements minimally affects the distribution of the eigenvalues.
\textcolor{black}{Therefore, if the diagonal elements of \textbf{M} follow a moderate-variance distribution with mean satisfying (\ref{eq:3.5}),  the ecological system will remain locally stable.} Nevertheless, when the variance  significantly surpasses that of the interspecific interactions, the  impact on stability depends on the specific patterns of the diagonal elements. When the self-regulation is stronger for species with fewer interactions, the impact of a considerable variance on stability remains negligible \cite{tang2014correlation}.

\subsubsection{Impacts of modularity}\label{sec:3.1.5}
The stability of  ecological systems can also be influenced by  network modularity. A modular network can be  divided into different modules or  subsystems, and the interactions within each subsystem are much more frequent than those between subsystems \cite{newman2006modularity}. 
The modularity of a network can be defined as
\begin{equation*}
    Q = \frac{L_w-\mathbb{E}[L_w]}{L_w+L_b},
\end{equation*}
where $L_w$ is the observed number of interactions within the subsystems, and $L_b$ is the number of inter-subsystem interactions.
Grilli et al.  studied the effect of modularity $Q$ on the stability of ecological systems \cite{grilli2016modularity}. \textcolor{black}{They found that modularity exhibits a moderate stabilizing effect when the subsystems have similar sizes and the overall mean interaction strength is negative. In particular, the stabilizing effect becomes stronger for negative correlations. Conversely, anti-modularity is highly destabilizing, except for the case where the overall mean interaction strength is close to zero. The authors further investigated the effect of modularity in food webs through numerical simulations and found that the results remain qualitatively unchanged \cite{grilli2016modularity}.
}

\subsubsection{Impacts of dispersal}\label{sec:3.1.6}
Spatial flows (e.g., exchanges of individuals, energy, and material) among local ecosystems are ubiquitous in nature \cite{polis1995extraordinarily}. Gravel et al. incorporated dispersal, i.e., spatial movement of species among local ecological systems, into the stability analysis of meta-ecological systems \cite{gravel2016stability}. The community matrix \textbf{M} of a meta-ecological system can be expressed as the sum of three matrices, i.e., 
\begin{equation*}
    \textbf{M} = \textbf{D} + \textbf{Q}+\textbf{K},
\end{equation*}
where $\textbf{D}$ is a  diagonal matrix that accounts for  intraspecific density dependence, $\textbf{Q}$ is a  matrix representing dispersal among patches with diffusion coefficient $q$, and $\textbf{K}$ is a block diagonal matrix that contains the local community matrices (random matrices with $S$ species, connectance $C$, and interspecific interaction strength $\sigma$). By assuming  the number of  local systems and $S$ are large, the authors obtained the following stability criterion \cite{gravel2016stability}:
\begin{equation}
    \sigma\sqrt{C(S-1)} < d+q.
\end{equation}
Therefore, the effect of dispersal can promote stability in meta-ecological systems, in which  dispersal can move the most of the eigenvalues of the community matrix towards more negative values and shrinks the range of the remaining in proportion to the number of effective patches. Interestingly, Baron et al. discovered that the introduction of dispersal can lead to the Turing-type instability (\textcolor{black}{i.e., the equilibrium point is unstable with respect to disturbances of a finite range of wavelengths}) in an ecological system with a trophic structure \cite{baron2020dispersal}. While the inclusion of trophic structures often enhances stability in large ecological systems (e.g., food webs) \cite{gross2009generalized,johnson2014trophic}, it can lead to a peak in the real part of the maximum eigenvalue of the community matrix in the presence of dispersal, rendering the equilibrium point unstable \cite{baron2020dispersal}.

\subsubsection{Impacts of time delay}\label{sec:3.1.7}
Empirical evidence has indicated that species interactions often exhibit time lags, rather than occurring instantaneously \cite{niebur1991collective,ohira1999resonance}. Thus, time delays can hold significant implications for stability and coexistence. Notably, Pigani et al. investigated delay effects on the stability of large ecological systems \cite{pigani2022delay},
which can be captured by the following linearized equation: 
\begin{equation}\label{eq: delay}
    \dot{\textbf{z}}(t) = \textbf{M}\textbf{z}(t) + \textbf{M}_{\text{delay}}\textbf{z}(t-\tau),
\end{equation}
where $\textbf{M}_{\text{delay}}$ is the community matrix with delay.  For simplicity, the community matrix $\textbf{M}$ is set to $-d_M \textbf{I}$ where $d_M\geq 0$, \textcolor{black}{which guarantees stability for a sufficiently small community matrix delay. In other words, the current intraspecific interactions are always stabilizing.} $\textbf{M}_{\text{delay}}$ is supposed to be a random matrix with a constant diagonal entry $-d_{\text{delay}}\leq 0$ and off-diagonal elements normally distributed with zero mean, standard derivation $\sigma_{\text{delay}}$, and connectance $C_{\text{delay}}$. The authors proved that 
(\ref{eq: delay}) is stable if all the roots of the characteristic equation, defined as $H(z)=z-\lambda - \lambda_{\text{delay}}e^{-z\tau}$ where $\lambda$ and $\lambda_{\text{delay}}$ are the eigenvalues of \textbf{M} and $\textbf{M}_{\text{delay}}$, respectively, have negative real parts \cite{pigani2022delay}. Importantly,  these roots can be solved implicitly as a function of $\tau$, along with the parameters $d_M$, $d_{\text{delay}}$, $\sigma_{\text{delay}}$, and $C_{\text{delay}}$. 
They found that an increasing delay tends to destabilize the system, and if a system is already unstable for $\tau=0$, the delay cannot stabilize the system. Furthermore, the authors determined the critical delay $\tau^*$ as the minimum value of $\tau$ above which the system becomes unstable. Finally, distributed delay was considered, i.e., the second term in (\ref{eq: delay}) is modified to $\textbf{M}_{\text{delay}}\int_{0}^{\infty} \exp{\{-\frac{\tau}{\tilde{\tau}}\}}/\tilde{\tau}\textbf{z}(t-\tau)d\tau $, \textcolor{black}{which can be derived from logistic or resource competition models \cite{gupta2021effective}.} They found that the system becomes more and more unstable as $\tilde{\tau}$ increases. However, if $\tilde{\tau}$ is large enough, the system eventually goes back to the stable regime \cite{pigani2022delay}.

\subsubsection{Limitations of model-implicit approaches}\label{sec:3.1.8}
May's approach and the follow-up studies  offer a valuable theoretical framework for understanding the stability of  ecological networks with various structures. Yet, these model-implicit approaches heavily rely on randomly sampling interactions without an underlying dynamical model, which results in a lack of  biologically realistic representation of species interactions. Consequently, the analytical results produced by these approaches are independent  from the underlying  model
from which they are hypothetically derived. While informative, these results may not accurately capture the nuanced and complex interactions that characterize the stability of real ecological systems. One way to enhance this framework is to integrate underlying ecological models into the linear stability analysis, as detailed in the next subsection.

\subsection{Model-explicit approaches}\label{sec:3.2}
\textcolor{black}{The most commonly used model to describe the dynamics of  ecological systems is the generalized Lotka-Volterra (GLV) model \cite{bhargava1989generalized,malcai2002theoretical} defined as
\begin{equation}\label{eq: glv}
    \dot{x}_i(t) = x_i(t)\Big[r_i  + \sum_{j=1}^S \textbf{A}_{ij} x_j(t)\Big]
\end{equation}
for $i=1,2,\dots, S$, where $r_i$ is the intrinsic growth rate of species $i$ and $\textbf{A}\in\mathbb{R}^{S\times S}$ is the interaction matrix whose off-diagonal elements $\textbf{A}_{ij}$ represent the effect that species $j$ has upon species $i$. By assuming the existence of a feasible equilibrium point ${\bf x}^*$ (i.e., $x_i^* \geq 0$ for all species), the community matrix can be computed as
$
    \textbf{M} = \textbf{X}\textbf{A},
$
where $\textbf{X}\in\mathbb{R}^{S\times S}$ the diagonal equilibrium matrix such that $\textbf{X}_{ii} = x_i^*$. Hence, the community matrix can be regarded as the scaled interaction matrix. Similarly to model-implicit approaches, we can investigate the stability properties of the GLV model based the linearized equation (\ref{eq: linearized}). 
}

\begin{figure}[t]
\centering
  \includegraphics[scale=0.65]{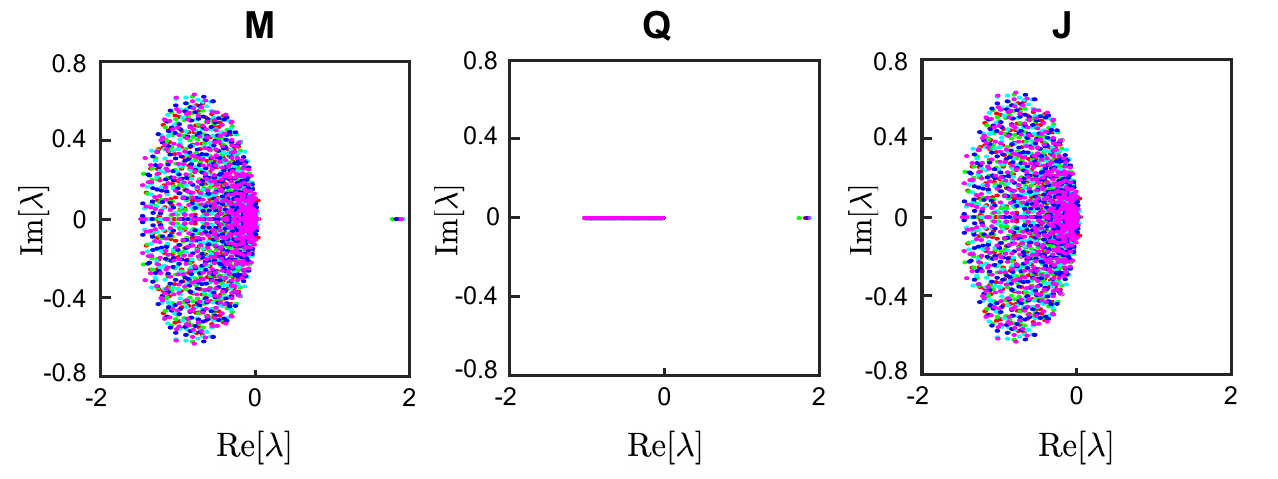}
  \caption{\textcolor{black}{Eigenvalue distributions of  the three matrices \textbf{M}, \textbf{Q}, and \textbf{J}. The off-diagonal elements of \textbf{A} are sampled from a normal bivariate distribution with identical marginal $\mu_A=5/S$, $\sigma_A=5/{\sqrt{S}}$, and correlation $\rho_A=-0.5$. The diagonal elements of \textbf{A} are fixed to $\mu_d = -1$, while the diagonal elements of   \textbf{X} are sampled from a uniform distribution on [0,1]. The two matrices \textbf{Q} and \textbf{J} are defined as $\textbf{Q}=\textbf{X}[(\mu _d-\mu_A)\textbf{I}+\mu_A \textbf{1}]$ and $\textbf{J}=\textbf{X}(\mu _d\textbf{I}+\textbf{B})$, respectively, such that the bulk of the eigenvalue distributions of \textbf{M} are the same as these of \textbf{J}, and the outlier eigenvalue of \textbf{M} is the same as that of  \textbf{Q}. Here, \textbf{B} is a random matrix with zero diagonal elements whose off-diagonal elements have mean zero and variance $\sigma_A^2$.  $S=500$. This figure was redrawn from  \cite{gibbs2018effect}.}}
  \label{fig: outlier}
\centering
\end{figure}

\subsubsection{Impacts of equilibria}\label{sec:3.2.1}
Equilibrium points play a significant role in the local stability analysis of the GLV model, as its community matrix is directly related to the equilibrium points of the system. Gibbs et al. explored the eigenvalue spectrum of the community matrix by assuming that the diagonal equilibrium matrix \textbf{X} and the interaction matrix \textbf{A} are randomly drawn from an arbitrary distribution with positive support and a bivariate distribution, respectively \cite{gibbs2018effect}. They proved that the eigenvalue spectrum of the community matrix \textbf{M} consists of a bulk of eigenvalues with mean $\mu_X(\mu_d-\mu_A)$ determined by the eigenvalues of the matrix $\textbf{J}=\textbf{X}(\mu _d\textbf{I}+\textbf{B})$ and an outlier eigenvalue determined by the largest eigenvalue of the matrix $\textbf{Q}=\textbf{X}[(\mu _d-\mu_A)\textbf{I}+\mu_A \textbf{1}]$ (\textbf{1} is the all-one matrix) such that
\begin{equation*}
    \lambda_{\text{outlier}} = \mu_X[\mu_d + (S-1)\mu_A],
\end{equation*}
where $\mu_X$, $\mu_d$, and $\mu_A$ represent the mean of the diagonal of \textbf{X}, the diagonal of \textbf{A}, and the off-diagonal of \textbf{A}, respectively, see Figure \ref{fig: outlier}. 
Based on their finding, the authors concluded that for mutualistic ecological systems, if the interaction matrix is stable, then the community matrix  will also be stable \cite{gibbs2018effect}. It is important to note that Gibbs et al. made the assumption that the distribution from which the equilibrium is drawn is independent of the interaction matrix, which does not typically hold in reality. Nevertheless,
Liu et al.  discovered that if the equilibrium point follows a specific distribution related to the elements of the interaction matrix \textbf{A}, Gibbs et al.'s assumption remains valid \cite{liu2023feasibility}.
Around the same time, Stone  demonstrated that for large ecological systems described by the GLV model, the stability of the interaction matrix implies the stability of the community matrix, under the condition that all species have positive equilibria \cite{stone2018feasibility}. 

\subsubsection{Impacts of extinction}\label{sec:3.2.2}
Extinction is generic in the GLV model with random interactions \cite{baron2023breakdown,emary2021can,pettersson2020stability}. Recent studies  have highlighted that an equilibrium point of the GLV model is locally stable if and only if all the eigenvalues of the reduced interaction matrix (i.e., the interaction matrix between the species in the surviving sub-community) have negative real parts \cite{barbier2021fingerprints,baron2023breakdown,biroli2018marginally,stone2018feasibility}. Baron et al.  focused on the eigenvalue spectrum of the reduced interaction matrix in the GLV model, where the spectrum of the reduced interaction matrix also consists of a bulk set of eigenvalues and an outlier \cite{baron2023breakdown}. Importantly, the authors demonstrated that the universality principle holds for the bulk region (i.e., following Ellipse Law) but not for the outlier eigenvalue. The outlier eigenvalue can be solved from the generating functional approach.
Prediction of extinction boundary and new equilibrium after a primary extinction were also discussed in \cite{pettersson2020stability} and \cite{emary2021can}, respectively. 

\begin{figure}[t]
    \centering
    \includegraphics[scale=1.5]{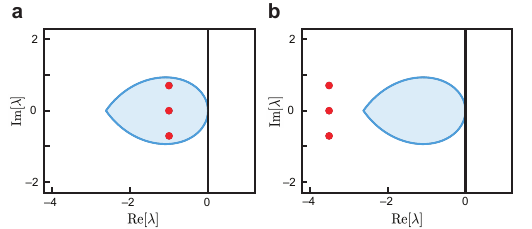}
    \caption{Eigenvalue distributions with time delay, where the blue teardrop-shaped area represents the stability region with time delay. (a)  The eigenvalues (represented by red dots) are within the teardrop-shaped area, indicating that the corresponding GLV model with time delay is locally stable. (b) The eigenvalues are outside the teardrop-shaped area, indicating that the corresponding GLV model with time delay is  unstable.  $\tau=1$. The figure was redrawn from \cite{yang2023time}.}
    \label{fig:tear}
\end{figure}

\subsubsection{Impacts of time delay}\label{sec:3.2.3}
Similar to Pigani et al.'s work as discussed in  Section \ref{sec:3.1.7}, Yang et al.  investigated the stability of the time-delayed GLV model \cite{yang2023time}, which is defined as
\begin{equation}\label{eq:tdglv}
    \dot{x}_i(t)=x_i(t)\Big[r_i+\sum_{j=1}^S\textbf{A}_{ij}x_j(t-\tau)\Big].
\end{equation}
The linearized equation of the time-delayed GLV model is given by
\begin{equation}\label{eq:dglv}
    \dot{\textbf{z}}(t) = \textbf{M}_{\text{delay}}\textbf{z}(t-\tau),
\end{equation}
where $\textbf{M}_{\text{delay}}=\textbf{X}\textbf{A}$ is the community matrix with delay. For simplicity, $\textbf{X}$ is set to $\textbf{I}$ such that each species has unit abundance. \textcolor{black}{However, the validation of this assumption remains an open problem \cite{yang2023time}.} They proved that (\ref{eq:dglv}) is stable if all the roots of the characteristic equation, defined as
$
    H(z) = z-\lambda_{\text{delay}} e^{-z\tau}
$
where $\lambda_{\text{delay}}$ are the eigenvalues of $\textbf{M}_{\text{delay}}$, have negative real parts. This equation implies that all the eigenvalues of $\textbf{M}_{\text{delay}}$ are required to be located in a teardrop-shaped region defined by $\text{Re}[z]<0$ to ensure stability, see Figure \ref{fig:tear}. Since the distribution of the eigenvalues of $\textbf{M}_{\text{delay}}$ (i.e., the interaction matrix \textbf{A} is this case) with random, mutualistic, or competitive interactions is well-established, it becomes straightforward to estimate $\max{(\text{Re}[z])}$. Based on the theoretical findings, the authors found that time delay plays a critical role in ecological systems, where large delay is often destabilizing, but short delay can considerably enhance community stability ($\tau=0.2$ and 0.5 were used for the small and large time delay in their simulations, respectively) \cite{yang2023time}. 

Similar conclusions were also obtained by Saeedian et al.  when investigating the impact of time delay on the emergent stability patterns within the time-delayed GLV model (\ref{eq:tdglv}) \cite{saeedian2022effect}. Significantly, they further determined the existence of a Hopf  bifurcation (\textcolor{black}{i.e., the two  complex  conjugate  eigenvalues  of  the community matrix, with non-zero imaginary part, simultaneously cross the imaginary axis into the right half-plane})  at $\tau_c$, which is computed as
\begin{equation*}
    \tau_c = \min_{\lambda_i\in \lambda(\textbf{M}_{\text{delay}})} \frac{1}{\lambda_i} \text{artan}\Bigg|\frac{\text{Re}[\lambda_i]}{\text{Im}[\lambda_i]}\Bigg|
\end{equation*}
for  the time-delayed  GLV model  with  a  large  number of  species.

\subsubsection{Impacts of stochastic noises}\label{sec:3.2.4}
\textcolor{black}{Real ecological systems are inherently stochastic with constant external perturbations and internal fluctuations \cite{tian2020estimating}.}
Krumbeck et al.  
developed a theoretical framework for the analysis of temporal stability of  ecological systems with stochastic noises \cite{krumbeck2021fluctuation}. They
provided analytical predictions for the power spectral density of the stochastic GLV model defined as
\begin{equation}\label{eq:sglv}
    \dot{x}_i(t) = x_i(t)\Big[r_i+\sum_{j=1}^S\textbf{A}_{ij}x_j(t)\Big] + \frac{1}{\sqrt{V}} \eta_i(t),
\end{equation}
where $V$ is the size of the species living domain, and $\eta_i(t)$ are Gaussian noises. 
The power spectral density is a statistical measure that can capture various aspects of temporal stability (e.g., the height of the spectrum gives information about the magnitude of stochastic fluctuations; the locations of nonzero peaks correspond to quasi-cyclic signals; a peak at zero indicates baseline wander). Specifically, the authors utilized the linearized equation of the stochastic GLV model (\ref{eq:sglv}), i.e.,
\begin{equation}\label{eq: s linearized}
    \dot{\textbf{z}}(t) = \textbf{M}\textbf{z}(t) + \boldsymbol{\zeta}(t),
\end{equation}
where $\boldsymbol{\zeta}(t)$ is a vector of Gaussian white noises with correlation matrix \textbf{B}. The power spectral density of fluctuations $\Psi(\omega)$ in the frequency domain can be computed as
\begin{equation*}
   \Psi(\omega) = (\textbf{M}-i\omega \textbf{I})^{-1}\textbf{B}(\textbf{M}^\top+i\omega\textbf{I})^{-1}, 
\end{equation*}
where $i$ here is the imaginary number. They  showed that different network structures (random,  mutualistic, and competitive) have unique signatures in the spectrum of fluctuations. They further investigated the effect of trophic structures and identified a gap in the power spectral density, which indicates that high-level trophic structures contribute to enhanced long-term temporal stability \cite{krumbeck2021fluctuation}.


\subsubsection{Impacts of evolved system size}\label{sec:3.2.5}
Previous work on the linear stability of the GLV model is concerned with ecological systems of predetermined and fixed sizes. Galla  explored the stability of  ecological systems with evolved system size using methods from statistical mechanics and the theory of disordered systems \cite{galla2018dynamically,martin1973statistical,mezard1987spin}. Specifically, he exploited generating functionals to derive the effective dynamics for the GLV model (with $r_i=1$) computed as
\begin{equation}
    \dot{x}(t) = x(t)\Big[1-x(t) + \mu_A M(t) + \rho_A\sigma_A^2 \int_0^t G(t, t{'})x(t{'}) dt{'}+\eta(t)\Big],
\end{equation}
where $M(t)$ is the average species concentration, $G(t, t{'})$ is the response function, and $\eta(t)$ is the Gaussian noise. As described previously, $\mu_A$, $\sigma_A$ and $\rho_A$ are the mean, standard derivation, and correlation of the off-diagonal elements of the interaction matrix \textbf{A}, respectively. The effective process characterizes the dynamics of a single representative species, denoted by $x(t)$, and captures the statistical behaviors of the ecological system. The linearized effective dynamics can be computed as
\begin{equation}
    \dot{z}(t) = x^{*}\Big[-z(t)+\rho_A\sigma_A^2\int_0^t G(t, t{'})z(t{'})dt{'} + v(t) + \zeta(t)\Big],
\end{equation}
where $z(t)$  denotes fluctuations about the equilibrium $x^{*}$, $v(t)$ is the  deviation of the noise in the effective process, and $\zeta(t)$ is the Gaussian white noise of unit amplitude. By performing Fourier transform with a focus on the long-time behavior of perturbations ($\omega=0$), the author established that the GLV model has stable equilibrium points in the limit of large population for 
    \begin{equation}
        \sigma_A < \sqrt{\frac{2}{(1+\rho_A)^2}}.
    \end{equation}
This result implies that predator-prey relationships enhance stability, while variability in species interactions promotes instability \cite{galla2018dynamically}, which aligns with prior findings as reported in \cite{tang2014correlation, yan2017degree}.

Poley et al.  applied Galla's approach to analyze the stability of the GLV model with hierarchical interactions \cite{poley2023generalized} defined as
\begin{equation}
    \dot{x}_i^a(t) = x_i^a(t)\Big[r_i^a+\sum_{b=1}^B\sum_{j=1}^{S^b}\textbf{A}_{ij}^{ab}x_j^b(t)\Big]
\end{equation}
for $i=1,2,\dots,S^a$, where $B$ is the total number of local systems, and $a,b=1,2,\dots, B$ denote the indices of the local systems with the associated system size $S^a$ and $S^b$. Their findings indicate that a strong hierarchical structure is stabilizing but reduces diversity within the effective dynamics of the GLV model \cite{poley2023generalized}. Similarly, Sidhom and Galla  adapted Galla's approach to study the GLV model with nonlinear feedback \cite{sidhom2020ecological} defined as
\begin{equation}\label{eq: glv feedback}
    \dot{x}_i(t) = x_i(t)\Big[r_i + g\big(\sum_{j=1}^S \textbf{A}_{ij}\,  x_j(t)\big)\Big]
\end{equation}
for $i=1,2,\dots,S$, where the function $g(u)=\frac{2au}{a+2|u|}$ denotes the nonlinear feedback ($a$ is the saturation parameter). This form of feedback was initially introduced to describe the predator's growth rate when interacting with prey. It is plausible that the benefits from extra prey will eventually saturate as prey numbers become large. The authors concluded that the stability and diversity of ecological systems improves with the introduction of nonlinear feedback \cite{sidhom2020ecological}.

\begin{figure}[t]
    \centering
    \includegraphics[scale=0.55]{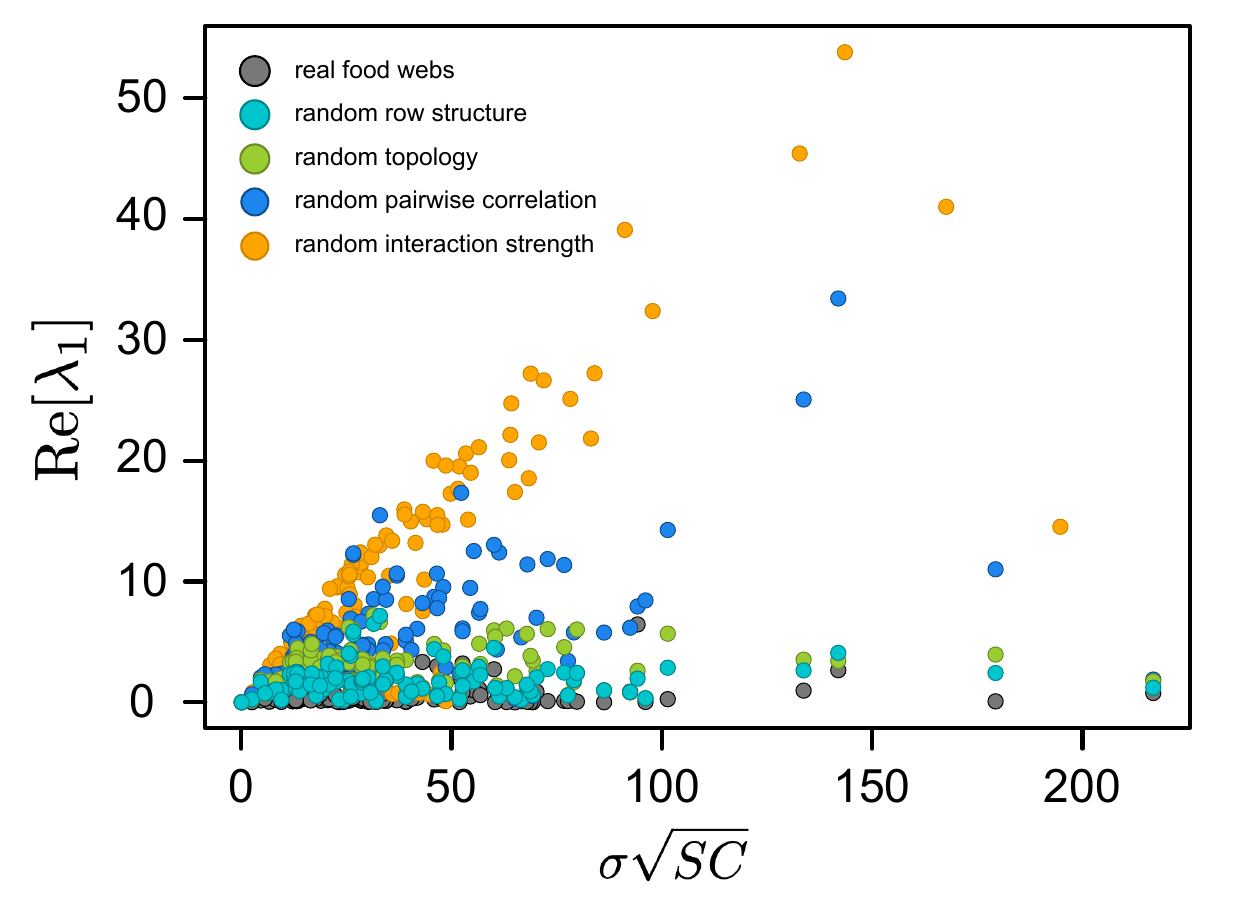}
    \caption{Stability of real and permuted food webs in relation to complexity. These permuted food webs were constructed from real food webs by removing some non-random features, namely the row structure, topology, pairwise correlation, and interaction strength. The more the ``rightmost''eigenvalues, i.e., $\text{Re}[\lambda_1]$, close to zero, the more stable of the food webs (because the diagonal elements of the community matrix are set to zero). This figure was redrawn from \cite{jacquet2016no}.}
    \label{fig:real}
\end{figure}

\subsubsection{Applications to real ecological systems}\label{sec:3.2.6}
Last but not least, while the stability analysis of random ecological systems offers powerful insights for the study of ecological dynamics and the understanding of how species interactions shape the stability and diversity of ecological communities, empirically 
applying the theory to real ecological systems poses a formidable challenge
\cite{ de1995energetics,emmerson2004predator,jacquet2016no,neutel2002stability,neutel2007reconciling,yodzis1981stability,yonatan2022complexity}. Notably, Jacquet et al.  performed the local stability analysis of 116 real food webs sampled worldwide from marine, freshwater, and terrestrial habitats using the GLV model \cite{jacquet2016no}. The community matrix \textbf{M} was constructed by multiplying the interaction matrix \textbf{A} with species biomass for each food web. \textcolor{black}{The interaction coefficients of \textbf{A} can be translated from the parameters of the corresponding Ecopath model, a trophic model that provides accurate representation of feeding interactions within food webs \cite{christensen1992ecopath}.} The authors then measured the stability of food webs using the real part of the dominant eigenvalue of the community matrix \textbf{M}. Using randomization tests, they found that negative correlation between interaction strengths with high frequency of weak interactions is a strong stabilizing property in real food webs, see Figure \ref{fig:real}. Their findings reveal that empirical food webs exhibit some non-random characteristics that lead to the absence of a complexity–stability relationship \cite{jacquet2016no}.

\section{Sign stability analysis}\label{sec:4}
\begin{figure}[t]
    \centering
    \includegraphics[scale=0.43]{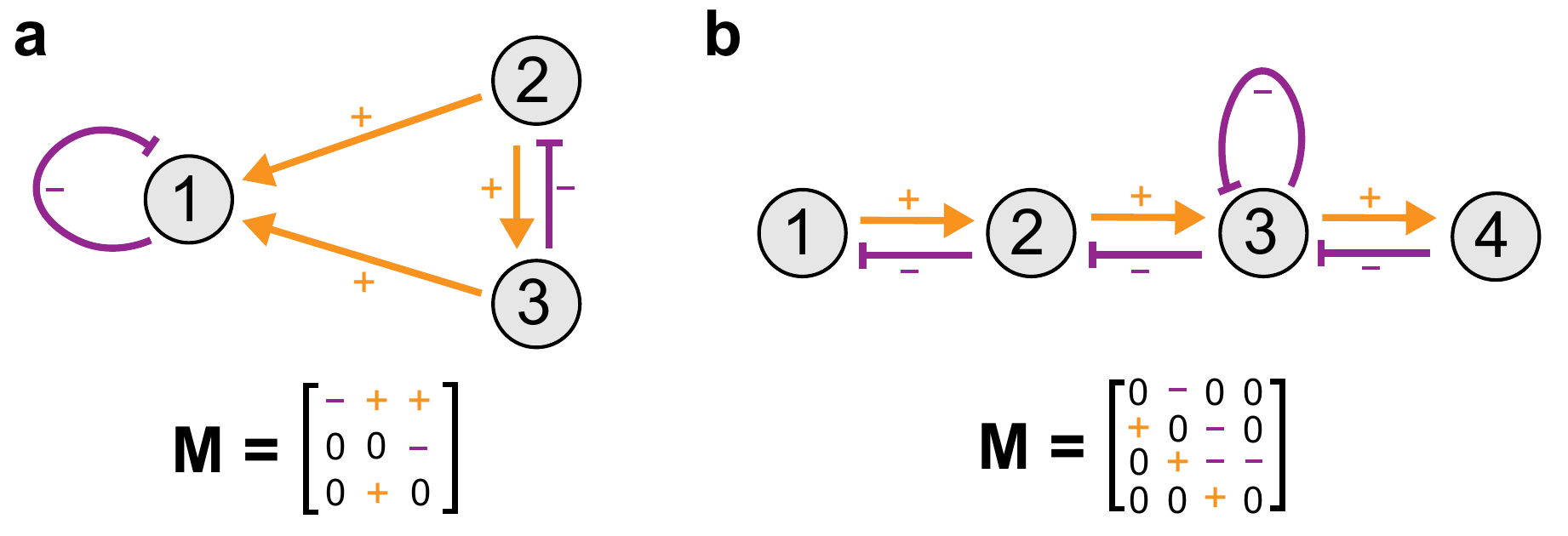}
    \caption{\textcolor{black}{Examples of signed ecological networks and their corresponding sign community matrices. (a) Self-regulated species 1 forms a commensal relationship with each species within the prey-predator pair of 2 and 3. (b) Linear trophic chain in which each successive species is preyed upon the subsequent one and species 3 is self-regulated. This figure was redrawn from \cite{logofet1982sign}.}}
    \label{fig:sign}
\end{figure}

The notion of sign stability is of particular interest for ecology, economics, chemistry, and engineering \cite{clarke1975theorems,eschenbach1998real,jeffries1977matrix,logofet1982sign,may1973qualitative,quirk1965qualitative,yedavalli2008ecological,yedavalli2009qualitative}. In an ecological system, the community matrix ${\bf M}$ (or the interaction matrix \textbf{A}) associated with certain models (e.g., the GLV model) might be only known qualitatively in the sense that the signs of the elements $\textbf{M}_{ij}$ can be determined with reasonable confidence, but the actual magnitudes may be very difficult to determine, see Figure \ref{fig:sign}. Such matrices are often referred to as sign matrices.  A sign matrix ${\bf M}$ is called sign stable (or sign semi-stable) if each of its eigenvalues has negative real part (or non-positive real part, respectively) for all numerical matrices of the same sign pattern \cite{Jeffries-74,jeffries1977matrix, Maybee-69,quirk1965qualitative,Yamada-Networks-90}.

\subsection{Characterizations of sign stable matrices}\label{sec:4.1}

A sign matrix ${\bf M} \in \mathbb{R}^{S\times S}$ is sign semi-stable if and only if the following conditions are all satisfied: (i) $\textbf{M}_{ii} \leq  0$ for
$i=1,2,\dots,S$; (ii) $\textbf{M}_{ij}\textbf{M}_{ji}\leq 0$ for $i,j=1,2,\dots, S$ and $i\neq j$; (iii) there exists no elementary cycle of length $k\ge 3$ in the digraph generated by \textbf{M} \cite{Maybee-69,quirk1965qualitative}. It is also known that if  $\textbf{M}_{ii} < 0$ for $i=1,2,\dots,S$, then (ii) and (iii) are necessary and sufficient for ${\bf M}$ to be sign stable. Quirk and Ruppert  further claimed that (i)-(iii) as well as (iv) $\textbf{M}_{ii}<0$ for some $i$ and (v) $\text{det}(\textbf{M})\neq 0$ are both necessary and sufficient for ${\bf M}$ to be sign stable \cite{quirk1965qualitative}. However, this was proved to be incorrect later by Jeffries et al. \cite{Jeffries-74}. 
Later on, Yamada  demonstrated that such an exceptional case is very rare and proved that the conditions (i)-(v) proposed by Quirk and Ruppert is actually necessary and sufficient for a system to be generically sign stable, i.e., sign stable for almost all parameter values except for some pathological cases with measure zero \cite{Yamada-Networks-90}. 

Quirk and Ruppert's conditions (i)-(v) can be translated into ecological terms to provide a comprehensive characterization of sign stable patterns \cite{logofet1982sign}. According to (i) and (iv), a sign stable ecological system should not include self-promoting species and must have at least one species with self-regulation. Condition (ii) necessitates the absence of both competitive and mutualistic relations. Condition (iii) states that there are no closed directed cycles with more than two edges in the network structure. Finally, condition (v) requires that  the digraph generated by \textbf{M}  must include a specific number of non-overlapping directed cycles that encompass all the nodes.

Jeffries et al.  introduced two additional conditions, distinct from (iv) and (v), called the coloring and matching conditions \cite{Jeffries-74}. The two conditions, combined with (i)-(iii), are necessary and sufficient for sign stability. Let $\mathcal{R}_M=\{i \text{ }| \text{ }\textbf{M}_{ii}\neq0 \}$. An $\mathcal{R}_M$-coloring of an undirected graph  is a partition of its nodes into two sets, black and white, such that each node in $\mathcal{R}_M$ is black (one of which can be empty), no black node has exactly one white neighbor, and each white node has at least one white neighbor. The coloring condition then states that in every $\mathcal{R}_M$-coloring of the undirected graph generated by \textbf{M}, all nodes are black. Additionally, a $(\mathcal{V}\sim \mathcal{R}_M)$-complete matching in an undirected graph  is a set $\mathcal{M}$ of  disjoint edges such that an exact cover of the node set can be obtained using the pairs in $\mathcal{M}$ and certain singletons from $\mathcal{R}_{\mathcal{M}}$. The matching condition then asserts that the undirected graph generated by \textbf{M} admits a $(V\sim\mathcal{R}_M)$-complete matching. However, the ecological systems characterized by Quirk and Ruppert \cite{quirk1965qualitative} or Jeffries et al. \cite{Jeffries-74} are not readily observable in nature.

\subsection{Applications to ecological systems}\label{sec:4.2}
The notion of sign stability has been applied to various ecological systems through different approaches. Dambacher et al.  introduced two qualitative metrics -- weighted feedback and weighted determinants based on the Hurwitz criterion \cite{dambacher2003qualitative}, which can be recast into two conditions: (i) the coefficients of the characteristic equation of \textbf{M} must have the same sign; (ii) the corresponding Hurwitz determinants must all be positive. The two metrics offer a practical mean to identify the relative degree to which stable parameter space can be constrained based on system structure and complexity. Remarkably, Haraldsson et al.  utilized the weighted feedback and weighted determinants to investigate the sign stability of social-ecological systems, which is an important tool to understand human-nature relations \cite{haraldsson2020model}. 


Moreover, Allesina and Pascual  studied the sign stability of random and empirical food webs through a fairly intuitive approach \cite{allesina2008network}. Given a community matrix \textbf{M} (either randomly generated or empirical), first determine the stability based on its eigenvalues. If it is stable (in terms of $\text{Re}[\lambda(\textbf{M})]\leq 0$), generate 100 matrices that have the same sign pattern with random magnitude. Lastly, measure the percentage of the random generated matrices that are stable. The percentage can tell whether the stability is due to a particular combination of coefficients or the sign pattern of the network itself. Based on this approach, the authors demonstrated that predator-prey interactions can promote stability, highly robust to perturbations of interaction strength, in real ecological systems \cite{allesina2008network}.

\section{Diagonal stability analysis}\label{sec:5}
The notion of diagonal stability was first introduced by Volterra in the 1930s \cite{volterra1959theory}. It has been  particularly useful for the stability analysis of  ecological systems and other networked systems \cite{ballantine1970stabilization,barker1978positive,berman1983matrix,bhaya1993discrete,kaszkurewicz1996diagonal,Kaszkurewicz-Book-00,pastravanu2010diagonal,sun2023gallery}. A matrix $\textbf{A}$ is called diagonally stable if there exists a diagonal matrix ${\bf P} \succ 0$
that renders 
\begin{equation*}
    \textbf{A}^\top \textbf{P} + \textbf{P} \textbf{A}= -\textbf{Q} \prec 0.
\end{equation*}
We use the notation $\mathcal{D}(\textbf{P})$  for the class of diagonally stable matrices. The positive diagonal matrix $\textbf{P}$ is often called Volterra multiplier in literature \cite{Redheffer-SIAM-85}. In many cases, the necessary and sufficient conditions for the Lyapunov stability of nonlinear systems are also the necessary and sufficient conditions for the diagonal stability of a certain matrix associated to the nonlinear system. This matrix naturally captures the underlying network structure of the nonlinear dynamical system. 

\textcolor{black}{
Diagonal stability has been successfully applied to various types of dynamical systems \cite{aleksandrov2018diagonal,khong2016diagonal,mason2006simultaneous,sootla2017block,sun2019diagonal,wu2021diagonal}.} It has nice ``structural consequences.'' For example, the principal sub-matrices of a diagonally stable matrix are also diagonally stable, implying that all the corresponding ``principal sub-systems'' of a given diagonally stable system are diagonally
stable \cite{Kaszkurewicz-Book-00}. \textcolor{black}{In this section, we begin by  introducing the theoretical framework for determining the Lyapunov stability of the GLV model via diagonal stability analysis. Subsequently, we delve into  the characterizations of diagonal stability of the interaction matrix that is associated with special interconnection or network structures, which offer an effective mean of determining the  stability of ecological systems.
}

\subsection{Lyapunov stability of the GLV model}\label{sec:5.1}
\textcolor{black}{
Before talking about the stability of the GLV model, we first introduce the notion of Persidskii-type systems.
}
Persidskii-type systems are typical examples that admit diagonal-type Lyapunov functions \cite{efimov2021analysis,efimov2019robust,Persidskii-ARC-69,platonov2023asymptotic}. We need to introduce a few concepts to define Persidskii-type systems.  A function $\textbf{f} : \mathbb{R}^S \to \mathbb{R}^S : \textbf{x} \mapsto \textbf{f}(\textbf{x})$ is called diagonal if, the $i$th component of $\textbf{f}$, i.e., $f_i$, is a function of $x_i$ alone. A function $f : \mathbb{R} \to \mathbb{R}$ is said to be in sector $[l, u]$ if $\forall x \in \mathbb{R}$, $f(x)$ lies between $lx$ and $ux$, i.e.,
\begin{equation*}
    (f(x)-ux)(f(x)-lx) \le 0.
\end{equation*}
For example, sector $[-1,1]$ means $|f(x)| \le |x|$. Sector $[0,\infty]$ means $f(x)$ and $x$ always have same sign. The class of infinite sector nonlinear functions $\mathcal{S}$ is defined to be functions
in sector $[0,\infty]$ that satisfy $\int_0^x f(\tau)  d\tau \to \infty$ as $|x| \to \infty$. Typical examples of infinite sector nonlinear functions are $f(x)=x$, $f(x)=x^3$, and $f(x)=\tanh(x)$ \cite{Kaszkurewicz-Book-00}.  A dynamic system $\dot{\textbf{x}}(t)=\textbf{f}(\textbf{x}(t))$ with $\textbf{x}(t) \in
\mathbb{R}^S$  and $\textbf{f}: \mathbb{R}^S \to \mathbb{R}^S$ is said to be of Persidskii-type if it has the following form:
\begin{equation}
\dot{x}_i(t) = \sum_{j=1}^S \textbf{A}_{ij} f_j\big(x_j(t)\big) \label{eq:Persidskii_ODE}
\end{equation}
for $i=1,2,\dots, S$, where $f_j \in \mathcal{S}$ for all $j=1,2,\dots, S$. In other words, $\textbf{f}$ is diagonal and $f_j$ is in the class of infinite sector nonlinear functions. Note that (\ref{eq:Persidskii_ODE}) can be used to describe a wide range of complex networked systems, where $\textbf{A}_{ij}$ capture the weighted wiring diagram. 

For Persidskii-type systems, we can introduce a diagonal-type Lyapunov function of the following form:
\begin{equation*}
V({\bf x}) = \frac 12 \sum_{i=1}^S p_i \int_0^{x_i} f_i (\tau)  d
\tau, \label{eq:Persidskii_V} 
\end{equation*}
where $p_i$ denotes the $i$th diagonal of \textbf{P}.
The equilibrium point $\textbf{x}^*={\bf 0}$ of the Persidskii-type system is globally asymptotically stable (\textcolor{black}{in the sense of Lyapunov}) if $\textbf{A} \in \mathcal{D}(\textbf{P})$.  This can be seen by computing $\dot{V}({\bf x})$
along the trajectory of (\ref{eq:Persidskii_ODE}), yielding 
\begin{equation*}
\dot{V}(\textbf{x}) = \textbf{f}^\top (\textbf{x}) ( \textbf{A}^\top \textbf{P} + \textbf{P} \textbf{A} ) \textbf{f}(\textbf{x}).
\end{equation*}
Since $\textbf{f}$ is diagonal and $f_j \in \mathcal{S}$, if $\textbf{A} \in\mathcal{D}(\textbf{P})$, therefore $\dot{V}(\textbf{x})$ is negative definite. Moreover, the functions $f_i  \in \mathcal{S}$ ensure the radial unboundedness of $V(\textbf{x})$. Hence, according to Lyapunov's theorems of stability, $\textbf{x}^*={\bf 0}$ is globally asymptotically stable (\textcolor{black}{in the sense of Lyapunov}).

\begin{figure}[t]
    \centering
    \includegraphics[scale=0.35]{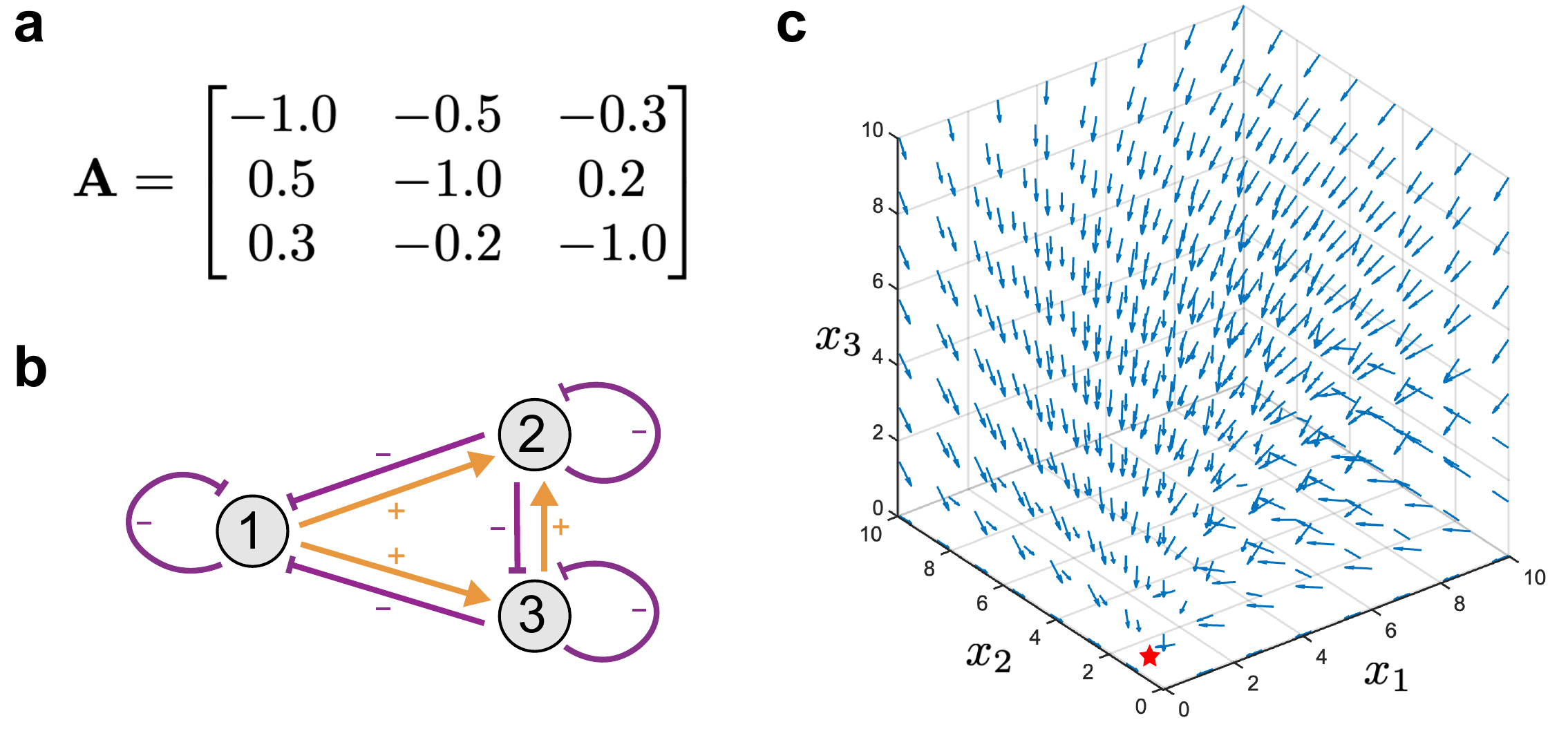}
    \caption{Example of a diagonally stable matrix. (a) Diagonally stable interaction matrix \textbf{A}. It is easy to check that the eigenvalues of $(\textbf{A}+\textbf{A}^\top)/2$ are all negative. 
(b) Corresponding ecological network associated with \textbf{A}. 
(c) Vector field plot of the corresponding GLV model with the non-trivial equilibrium point indicated by the red star. }
    \label{fig:diagonal}
\end{figure}

By assuming the existence of a  non-trivial equilibrium point ${\bf x}^*$
(i.e., $x_i^* >0$ for all species) and defining 
\begin{equation*}
    z_i(t) = \log\frac{x_i(t)}{x_i^*} \quad \text{and}\quad g_j\big(z_j(t)\big) = x_j^* (e^{z_j(t)}-1),
\end{equation*}
we can bring the GLV model (\ref{eq: glv}) into the form of Persidskii-type dynamics: 
\begin{equation*}\label{eq:GLV_z}
    \dot{z}_i(t) = \sum_{j=1}^S \textbf{A}_{ij} \, g_j \big(z_j(t)\big),
\end{equation*}
which admits the following diagonal-type Lyapunov function: 
\begin{equation*}
    V({\bf z}) = \sum_{i=1}^S p_i \int_0^{z_i} g_i (\tau) d \tau \quad\text{or}\quad V({\bf x}) = \sum_{i=1}^S p_i \left[ x_i - x_i^* - x_i^* \log \frac{x_i}{x_i^*} \right].
\end{equation*}
Using the above diagonal-type Lyapunov function, Goh  first proved that for the GLV model (\ref{eq: glv}), if the interaction matrix $\textbf{A}$ is diagonally stable, then the non-trivial equilibrium point $\textbf{x}^*$ in the positive orthant is globally asymptotically stable (in the sense of Lyapunov) \cite{goh1977global}. 
Therefore, diagonal stability allows the GLV model to lie in the unique fixed-point phase, as depicted by Bunin \cite{bunin2017ecological}. Here, we provide an example of a 3-by-3 diagonally stable matrix and the corresponding vector field plot of the GLV model, see Figure \ref{fig:diagonal}. Clearly, the non-trivial equilibrium point in the positive orthant is globally asymptotically stable.

The diagonal stability analysis has been applied to analyze the stability of various GLV-derived models. Following  Goh's findings, Wörz-Busekros  offered a sufficient condition for the global stability of the GLV model with continuous-time delay via diagonal stability analysis \cite{worz1978global}. Later on, Beretta and Takeuchi   considered the global asymptotic stability of diffusion models with multiple species heterogeneous patches, in which each patch is governed by the GLV model with continuous-time delay \cite{beretta1988global}. Concurrently, Beretta and Solimano  generalized the Wörz-Busekros's outcome by considering a non-negative linear vector function of the species \cite{beretta1988generalization}.
In addition, Kon  exploited the diagonal stability theory to determine the stability of the GLV model with an age structure \cite{kon2011age}. The notion of diagonal stability was further applied to more general ecological models such as Kolmogorov systems \cite{hou2015global}, where the GLV model is encompassed as a specific instance.  Significantly, the stability properties of many quasi-polynomial dynamical systems, often used to represent many biochemical processes, can be studied  through an equivalent GLV model that has a much simpler form \cite{hernandez1997lotka,motee2012stability}. This is based on the fact that the former can be transformed into the latter with some appropriate changes of variables.

\subsection{Characterizations of diagonally stable matrices}\label{sec:5.2}

A general characterization of the diagonally stable interaction matrix in the GLV model remains elusive for more than three species \cite{cross1978three,logofet2005stronger}, though there exist efficient optimization-based algorithms to numerically check if a given matrix is diagonally stable, e.g., polynomial-time interior point algorithms \cite{Boyd-Book-1994}. In cases where the dimension is three or less, the diagonal stability of \textbf{A} can be determined by examining the signs of its principal minors \cite{cross1978three}. For high-dimensional matrices under very special structural assumptions, one can derive necessary and sufficient conditions for diagonal stability. For example, if $\textbf{A}$ is Metzler (i.e., $\textbf{A}_{ij} \ge 0 \text{ for } i \neq j$), then $\textbf{A}$ is diagonally stable if and only if all principal minors of $-\textbf{A}$ are positive \cite{Berman-Book-94}.  More special examples are discussed in \cite{Redheffer-SIAM-85}. Moreover, there are approaches which can reduce the problem of determining whether $\textbf{A}\in\mathbb{R}^{S\times S}$ is diagonally stable into two simultaneous problems of  $(S-1)\times (S-1)$ matrices, but the method becomes intractable for large $S$ \cite{Redheffer-SIAM-85}.

Grilli et al.  imposed a condition of negative definiteness on \textbf{A} (i.e., all the eigenvalues of $\textbf{A}+\textbf{A}^\top$ are negative) to ensure the diagonal stability of \textbf{A} for studying the feasibility and coexistence of large random ecological systems \cite{grilli2017feasibility}. It is known that a negative definite matrix is also diagonally stable, and the condition is much easier to verify and characterized for random matrices. Later on, Gibbs et al. \cite{gibbs2018effect} estimated the ``rightmost'' eigenvalue of $(\textbf{A}+\textbf{A}^\top)/2$ with random interactions and proved that if 
\begin{equation}
    \mu_d+\sqrt{2S\sigma_A^2(1+\rho_A)}<0, 
\end{equation}
then \textbf{A} is diagonally stable, where $\mu_d$ is the mean of the diagonal elements of \textbf{A}, and  $\sigma_A^2$ and $\rho_A$ are the variance and correlation of the off-diagonal elements of \textbf{A}.

Recently, necessary and sufficient diagonal stability conditions for matrices associated with special interconnection or network structures were  studied \cite{arcak2011diagonal,arcak2008passivity,arcak2006diagonal, 4434115,simpson2021diagonal,wang2012diagonal}. If an ecological system described by the GLV model exhibits these  network structures, it will be effective and efficient to determine its Lyapunov stability through the diagonal stability of the interaction matrix \textbf{A}. We review these network structures as follows.

\begin{figure}[t]
    \centering
    \includegraphics[scale=0.52]{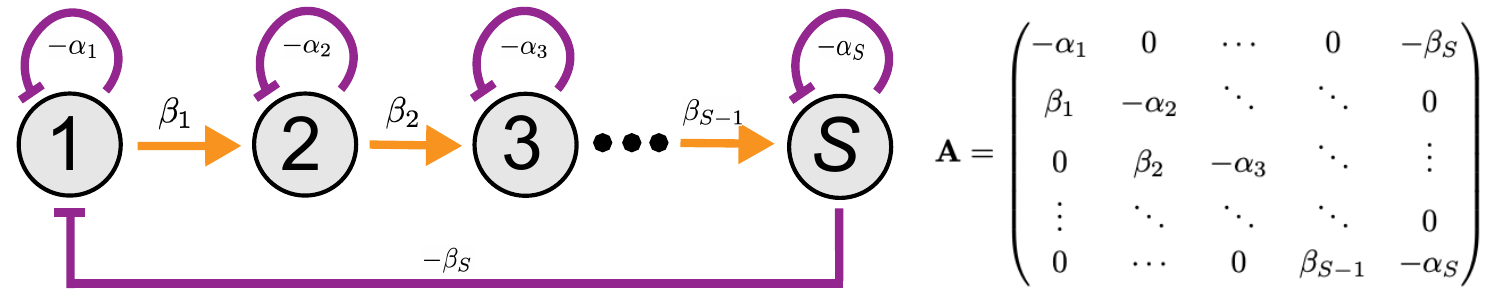}
    \caption{Negative feedback cyclic structure with the corresponding interaction matrix \textbf{A}, where $\alpha_i >0$, $\beta_i >0$ for $i=1, 2,\dots, S$. }
    \label{fig:negcircle}
\end{figure}

\subsubsection{Negative feedback cyclic structure}\label{sec:5.2.1}
In the negative feedback cyclic structure, the intermediate species play a facilitating role for the subsequent species, while the final species exerts inhibitory effects on the initial species, see Figure \ref{fig:negcircle}. It has been shown that the interaction matrix $\textbf{A}$ is Hurwitz, i.e., $\text{Re}[\lambda(\textbf{A})] < 0$, if it satisfies the so-called secant criterion \cite{Thron-BMB-91,Tyson-Book-78,wimmer2009diagonal}, i.e.,
\begin{equation}
    \frac{\beta_1\beta_2 \cdots \beta_S}{\alpha_1 \alpha_2\cdots \alpha_S} < 
\left[\sec\left(\frac{\pi}{S}\right)\right]^S. \label{eq:Secant}
\end{equation}
When $\alpha_i$ are equal, the  secant criterion (\ref{eq:Secant}) is also necessary for $\textbf{A}$ to be Hurwitz. Surprisingly, this secant criterion derived for linear stability is also a necessary and sufficient condition for diagonal stability of the corresponding class of matrices \cite{arcak2006diagonal}.

Note that a simple necessary condition for  $\textbf{A}$ to be diagonal stable is that all the diagonal elements $\textbf{A}_{ii}$ be negative. Since scaling the rows of $\textbf{A}$ by positive constants does not change its diagonal stability, one can safely assume $\textbf{A}_{ii}=-1$. Moreover, a reducible matrix $\textbf{A}$ can always be transformed into an upper-triangle form with a suitable permutation, and permutation does not change diagonal stability. Considering these two points, the secant criterion can be further generalized as follows \cite{arcak2011diagonal}.  Any matrix $\textbf{A}$ that can be transformed via a suitable permutation $\textbf{P}$ to the form of 
\begin{equation*}
    \textbf{P} \textbf{A} \textbf{P}^\top=\begin{pmatrix}
-1 & 0 & \cdots & 0 & -\tilde{\beta}_S \\
\tilde{\beta}_1   & -1 & \ddots & \ddots & 0 \\
0 & \tilde{\beta}_2   & -1 & \ddots & \vdots \\
\vdots & \ddots & \ddots & \ddots & 0 \\
0 & \cdots & 0 & \tilde{\beta}_{S-1} & -1
\end{pmatrix}
\end{equation*}
is diagonally stable if and only if 
$
|\gamma| \Phi(\text{sgn}(\gamma), S) < 1
$
where $\gamma=\tilde{\beta}_1\tilde{\beta}_2\cdots\tilde{\beta}_S\neq 0$ is the cycle gain and $\Phi(\text{sgn}(\gamma),S)=\cos^{S}(\pi/S)$ for $\gamma<0$; $1$ for $\gamma>0$. 

\subsubsection{Cactus structure}\label{sec:5.2.2}
Arcak further generalized the secant criterion to multiple cycles when the weighted digraph of the interaction matrix $\textbf{A}$, denoted as $G(\textbf{A})$, possesses a ``cactus'' structure, i.e., any pair of distinct simple cycles have at most one common node \cite{arcak2011diagonal}.  
Apparently, the negative feedback cyclic structure corresponds to a single cycle and is just a special case of the cactus structure. The digraph $G(\textbf{A})$ is defined to represent the off-diagonal entries of $\textbf{A}$ that has $S$ nodes and there is a directed edge $(j \to i)$ with weight $\textbf{A}_{ij}$ if and only if $\textbf{A}_{ij}\neq 0$. Self-loops corresponding to the diagonal entries $\textbf{A}_{ii}$ are excluded from $G(\textbf{A})$. 

\begin{figure}[t!]
\centering
\includegraphics[scale=0.34]{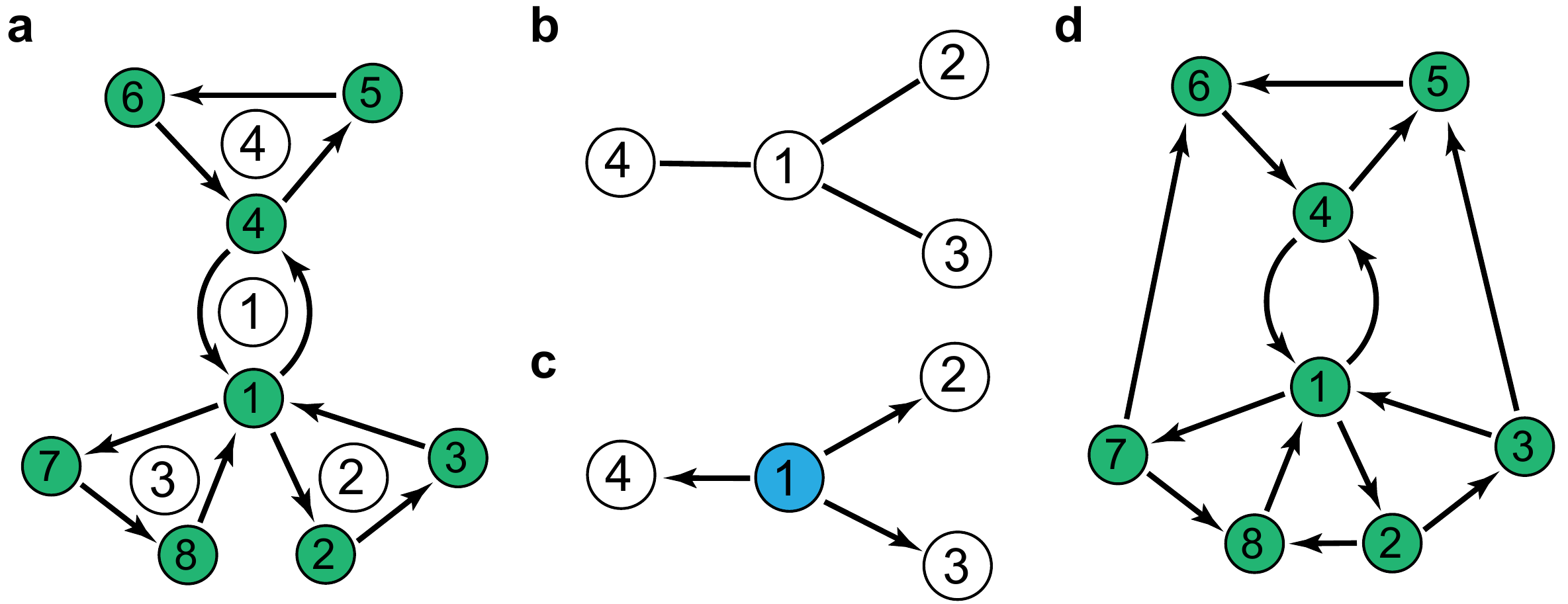}
\caption{Cactus and connected circle structures. (a) Digraph $G(\textbf{A})$ corresponding to the matrix $\textbf{A}$ with simple cycles $\{1,4,1\}, \{1,2,3,1\},\{1,7,8,1\}$, and $\{4,5,6,4\}$ labeled as 1, 2, 3 and 4, respectively. Note that any pair of the cycles have at most one common node. (b) Undirected graph  that describes which cycles in $G(\textbf{A})$ intersect. (c) Arborescence  constructed on the undirected graph in (b) according to the broadcasting algorithm with cycle 1 selected as the root (shown in cyan). (d) Digraph with the connected circle structure where any pair of the cycles have at most one common edge or one common node. Panels a, b, c were redrawn from \cite{arcak2011diagonal}.
}
\label{fig:cactus-Murat}
\end{figure}

Arcak made two assumptions about $\textbf{A}$ without loss of generality: (i)  $\textbf{A}_{ii}=-1$ for $i=1,\dots, S$; (ii) $G(\textbf{A})$ is strongly connected, or equivalently, $\textbf{A}$ is irreducible. Arcak then defined an undirected graph for describing which cycles of the digraph $G(\textbf{A})$ intersect. Subsequently, he constructed a spanning tree in the undirected graph and simultaneously assigned directions to its edges to create an arborescence, i.e., a digraph in which every node can be reached from the root by one and only one path, see Figure \ref{fig:cactus-Murat}a, b, c. Using the hierarchy of the cycles of $G(\textbf{A})$ established by this arborescence presentation, he sequentially generated dinequalities for the gains of each simple cycle.

Consider there are $l$ simple cycles in $G(\textbf{A})$. The length of cycle $j$ is denoted as $n_j$ for $j=1,2,\dots, l$. Define $\mathcal{I}_j = \{i_j^{(1)}, i_j^{(2)}, \dots, i_j^{(n_j)}\}$ and $J_i=\{1\le j \le l \text{ }|\text{ } i \in \mathcal{I}_j\}$ as the set of nodes traversed by cycle $j$ and the set of cycles that node $i$ belongs to, respectively. Denote $\gamma_j$ as the gain for cycle $j$. Then the stability condition can be summarized as follows: the interaction matrix $\textbf{A}$ satisfying the above two assumptions (i-ii) is diagonally stable if and only if there exist constants $\theta_{j}^{(i)} >0$ such that 
\begin{equation}
\begin{cases}
    \prod_{i \in \mathcal{I}_j} \theta_{j}^{(i)} > |\gamma_j|\Phi(\text{sgn}(\gamma_j), n_j) \text{ for } j=1,2, \dots, l \\
\sum_{j \in \mathcal{J}_i} \theta_{j}^{(i)} = 1 \text{ for } i=1,2, \dots, S
\end{cases}.
\end{equation}
\textcolor{black}{
Arcak then outlined a systematic procedure for constructing Lyapunov functions based on the above  stability condition. Notably, he illustrated the procedure with the GLV model for ecological systems \cite{arcak2011diagonal}. 
}

\textcolor{black}{
Later on, Wang and Nešić  extended the small gain condition to a more general circle structure called connected circles, where each pair of distinct simple cycles have at most one common edge or a common node \cite{wang2012diagonal}, see Figure \ref{fig:cactus-Murat}d. In addition to the two assumptions (i) and (ii) made by Arcak, the authors further assumed that (iii) \textbf{A} has non-negative off-diagonal elements. Define $\mathcal{L}_j=\{e_j^{(1)}, e_j^{(2)}, \dots, e_j^{(n_j-1)}\}$ as the set of edges traversed by cycle $j$. Then the stability condition can be modified as follows: the interaction matrix \textbf{A} satisfying the above three assumptions (i-iii) is diagonally stable if and only if there exist constants $\theta_{j}^{(i)} >0$  such that 
\begin{equation}
\begin{cases}
    \prod_{i \in \mathcal{I}_j} \theta_{j}^{(i)} > \gamma_j \text{ for } j=1,2, \dots, l \\
\sum_{j \in \mathcal{J}_i} \theta_{j}^{(i)} = 1 \text{ for } i=1,2, \dots, S
\end{cases},
\end{equation}
where $\gamma_j=\prod_{(p,q)\in \mathcal{L}_j}\textbf{A}_{pq}$ is the gain for cycle $j$. 
}

\subsubsection{Rank-one structure}\label{sec:5.2.3}
Recently, Simpson-Porco and Monshizadeh   explored necessary and sufficient conditions for the diagonal stability of \textbf{A} with a rank-one network structure \cite{simpson2021diagonal} defined as 
\begin{equation*}\label{eq:rk-one}
    \textbf{A} = -\textbf{D} + \textbf{S},
\end{equation*}
where $\textbf{D}=\text{diag}(d_1,d_2,\dots,d_S)$ is a positive diagonal matrix, and $\textbf{S}=\textbf{x}\textbf{y}^\top$ is a rank-one matrix for some \textbf{x}, $\textbf{y}\in\mathbb{R}^S$. Suppose that $\textbf{y}$ is a non-negative vector. The authors proved that the interaction matrix \textbf{A} with rank-one structure is diagonally stable if and only if 
\begin{equation}
    \sum_{i=1}^S\frac{1}{d_i}[x_iy_i]_{+}<1,
\end{equation}
where $x_i$ and $y_i$ are the $i$th elements of \textbf{x} and \textbf{y}, respectively, and $[\cdot]_{+}=\max{(\cdot, 0)}$ is the maximum operator with respect to zero. Significantly, they provided a theoretical stability analysis of automatic generation control in an interconnected nonlinear power system based on the above condition. While the rank-one structure may not naturally occur in real ecological networks, it might be useful in synthetic ecological networks or specific ecological scenarios.

\section{D-stability analysis}\label{sec:6}
\begin{figure}[t]
    \centering
    \includegraphics[scale=0.5]{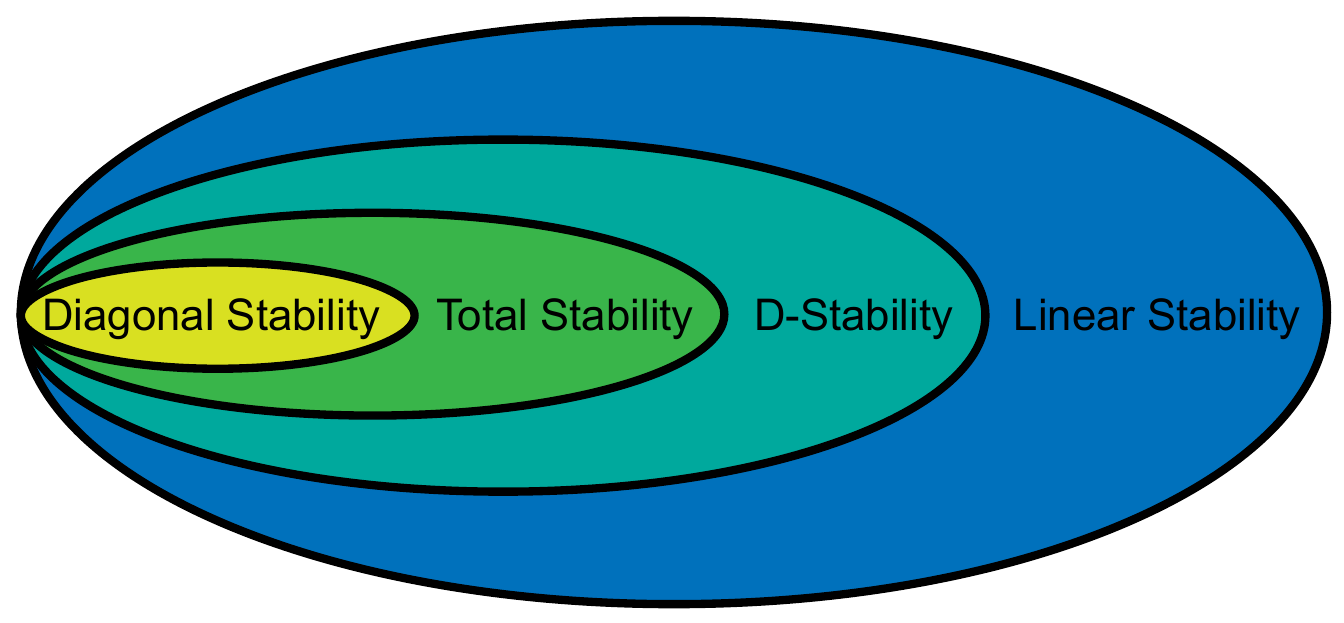}
    \caption{Schematic containment plot of matrix stability notions. Detailed proofs of the relations among these stability notions can be found in \cite{logofet2005stronger}. This figure was redrawn from \cite{saeedian2022effect}.}
    \label{fig:total}
\end{figure}

The notion of D-stability was originally introduced by Arrow and McManus \cite{arrow1958note} and Enthoven and Arrow \cite{enthoven1956theorem} in the late 1950s. D-stability can be considered as a weaker version of diagonal stability, see Figure \ref{fig:total}. A matrix \textbf{A} is called D-stable if for any positive diagonal matrix \textbf{X}, the matrix \textbf{X}\textbf{A} is stable (\textcolor{black}{in terms of $\text{Re}[\lambda(\textbf{X}\textbf{A})]\leq 0$)}. Clearly, the definition of D-stability is particularly relevant to the GLV model, where its community matrix is represented as $\textbf{M}=\textbf{X}\textbf{A}$. In this section, we first investigate the characterizations of D-stable matrices and then discuss existing results regarding D-stability in the context of the GLV model.

\subsection{Characterizations of D-stable matrices}\label{sec:6.1}
Explicitly characterizing D-stable matrices is only known for dimension less than or equal to four \cite{logofet2005stronger}. For example, a matrix $-\textbf{A}\in\mathbb{R}^{4\times 4}$ is D-stable if and only if all the principal minors of \textbf{A} are positive and $f(x,y,z)>0$ for all positive $x$, $y$, and $z$, where 
\begin{equation*}
    f(x, y, z) = E_1(\textbf{D}\textbf{A})E_2(\textbf{D}\textbf{A})E_3(\textbf{D}\textbf{A})-E_3(\textbf{D}\textbf{A})^2-\textbf{E}_1(\textbf{D}\textbf{A})^2E_4(\textbf{D}\textbf{A}),
\end{equation*}
and $\textbf{D}=\text{diag}(1,x,y,z)$. Here, $E_i(\cdot)$ denotes the sums of the order $i$ principal minors of a given matrix. However, the problem becomes highly challenging as the dimension grows larger \cite{bhaya1993discrete,gibbs2018effect,logofet2005stronger,Redheffer-SIAM-85}. Nonetheless, numerous sufficient conditions have been proposed \cite{johnson1974sufficient,leite2003improved,peaucelle2000new}. Some of the better-known ones are:
\begin{itemize}
    \item If \textbf{A} is diagonally stable,  then it is D-stable \cite{arrow1958note}. In fact, it is essentially the condition that Arrow and McManus offered.
    \item If \textbf{A} is Metzler (i.e., $\textbf{A}_{ij}\geq 0$ for $i\neq j$) and all the principal minors of $-\textbf{A}$ are positive,  then it is D-stable \cite{fiedler1962matrices}.
    \item If there exists a positive diagonal matrix \textbf{D} such that $-\textbf{A}\textbf{D} = \textbf{B}$ satisfies $\textbf{B}_{ii} > \sum_{j\neq i} |\textbf{B}_{ij}|$ for $i=1,2,\dots,S$, then \textbf{A} is D-stable \cite{quirk1968introduction}. The matrix $-\textbf{A}$ is referred to as quasi-dominant diagonal. 
    \item If \textbf{A} is triangular with $\textbf{A}_{ii}<0$ for $i=1,2,\dots,n$, then it is D-stable \cite{johnson1974sufficient}.  This is the most straightforward condition for D-stability. 
    \item If \textbf{A} is sign stable, then it is D-stable \cite{johnson1974sufficient}. 
    \item The element-wise product of $\textbf{P}$ and \textbf{A} is stable for each positive definite symmetric matrix \textbf{P} \cite{r1974hadamard}.
\end{itemize}
Regrettably,  none of these conditions is necessary for D-stability. More sufficient conditions can be found in \cite{johnson1974sufficient}.


\subsection{Results with the GLV model}\label{sec:6.2}
Here, we present a series of findings and speculation regarding D-stability in the context of analyzing the GLV model.
Chu  offered a solvable Lie algebraic condition for the equivalence of four stability notions -- linear stability, D-stability, total stability (defined in  Section \ref{sec:7}), and diagonal stability for the GLV model \cite{chu2007equivalent}. Given a Lie algebra $\mathbb{L}$, let's define the inductive sequence as 
\begin{equation*}
    \mathbb{L}^{(0)}=\mathbb{L}, \quad \mathbb{L}^{(i+1)}=\{\textbf{A}\textbf{B}-\textbf{B}\textbf{A} \text{ } |\text{ } \textbf{A}, \textbf{B}\in\mathbb{L}^{(i)}\}.
\end{equation*}
The Lie algebra $\mathbb{L}$ is referred to as solvable if there exists a positive integer $k>0$ such that $\mathbb{L}^{(k)}=\{\textbf{0}\}$. The author proved that given the GLV model (\ref{eq: glv}), if the two matrices, defined as $(\textbf{A}+\textbf{A}^\top)/2$ and $(\textbf{A}-\textbf{A}^\top)/2$, generate a solvable Lie algebra, then the linear stability, D-stability, total stability, and diagonal stability of the interaction matrix \textbf{A} are equivalent. Moreover, although diagonal stability is not equivalent to Lyapunov stability in general, based on intensive numerical simulations, Rohr et al. conjectured that for mutualistic ecological systems captured by the GLV model, if the interaction matrix \textbf{A} is stable (in the sense of Lyapunov), then \textbf{A} is D-stable \cite{rohr2014structural}. This conjecture has an important consequence for modeling mutualistic ecological systems with the GLV model, implying that if \textbf{A} is  stable (in the sense of Lyapunov), any feasible equilibrium point is globally stable, corresponding to the unique fixed-point phase, as illustrated by Bunin \cite{bunin2017ecological}. Similarly, Gibbs et al. observed that large random matrices are D-stable almost surely (i.e., the set of positive diagonal matrices leading to instability has measure zero), so feasible unstable equilibrium points are very unlikely for the GLV model with random interactions \cite{gibbs2018effect}.

\section{Total Stability}\label{sec:7}
The notion of D-stability is also closely linked to another stability notion termed total stability.  A matrix \textbf{A} is called totally stable if any principal sub-matrix of \textbf{A} is D-stable \cite{hofbauer1998evolutionary,quirk1965qualitative}. Any principal sub-matrix of a totally stable matrix is also totally stable. In addition, total stability is a stronger notion compared to D-stability, but is weaker than diagonal stability, as illustrated in Figure \ref{fig:total}. Consequently, any totally stable matrix is inherently D-stable, and the class of totally stable matrices is closed under transposition and multiplication by a positive diagonal matrix. This implies that total stability holds simultaneously for both the interaction matrix \textbf{A} and the community matrix \textbf{M} of the GLV model \cite{logofet2018matrices}. A direct way to characterize a totally stable matrix is to check the D-stability of all its principal sub-matrices, which is computationally expensive. However, explicit characterizations are only known for $S\leq 3$ \cite{logofet2018matrices,logofet1987on}, similar to diagonal stability and D-stability. A necessary condition of total stability is that all the principal minors of \textbf{A} of odd orders are negative, while those of even orders are positive \cite{logofet2018matrices}. Here, we provide an example of a 3-by-3 totally stable matrix that is not diagonally stable as follows (borrowed from \cite{logofet2018matrices}):
\begin{equation*}
    \textbf{A} = \begin{bmatrix}
        -1 & -0.3 & -2\\
        -0.95 & -1 & -1.92\\
        -0.49 & -0.495 & -1
    \end{bmatrix}.
\end{equation*}
It can be shown that all principal sub-matrices of \textbf{A} are D-stable. For example, the real parts of the eigenvalues of the following  product
\begin{equation*}
    \textbf{X}\textbf{A}=\begin{bmatrix}
        x_1^* & 0\\0 & x_2^*
    \end{bmatrix}
    \begin{bmatrix}
        -1 & -0.3\\
        -0.95 & -1
    \end{bmatrix}
\end{equation*}
are equal to $-x_1^*/2-x_2^*/2$. On the other hand, $(\textbf{A}+\textbf{A}^\top)/2$ contains a positive eigenvalue, so \textbf{A} is not diagonally stable.

\textcolor{black}{Ecologically, the principal sub-matrix of the interaction matrix  refers to the interactions between the subset of species that remains after removal or extinction of certain species. Therefore, total stability is related to the so-called ``species-deletion stability'' concept introduced in \cite{pimm1982food}, implying that the property is preserved after the removal or extinction of any group of species from the initial composition \cite{logofet2005stronger}. In the context of the GLV model, the total stability of the interaction matrix suggests that the remaining species will reach a new locally stable equilibrium after the removal or extinction of certain species.}

\section{Sector stability analysis}\label{sec:8}
\textcolor{black}{
The notion of sector stability, proposed by Goh, focuses on the stability properties of semi-feasible equilibrium points (i.e., $x_i^*=0$ for some $i$) for studying  interesting processes in ecology, e.g., succession and extinction  \cite{case1979global,goh1978sector}. A semi-feasible equilibrium point   is called (locally) sector stable if every trajectory of the system that starts within a non-negative neighborhood remains in the same or an even larger non-negative neighborhood and eventually converges to that equilibrium point.
The definition  is analogous to that of (local) asymptotic stability. However, sector stability restricts the trajectories within a non-negative part of an open neighborhood of the equilibrium point.  
}

\textcolor{black}{
Let's consider a generic population dynamics model defined as
\begin{equation}\label{eq:gpm}
    \dot{x}_i=x_if_i(x_1,x_2,\dots,x_S)
\end{equation}
for $i=1,2,\dots, S$, where $f_i$ have continuous partial derivatives at every finite point in the state space. Let $I_0=\{i\text{ }|\text{ }x_i^*=0\}$ and $I_1=\{i\text{ }|\text{ }x_i^*>0\}$. Goh proved the semi-feasible equilibrium point $\textbf{x}^*$ of the generic population model (\ref{eq:gpm}) is locally sector stable if all the eigenvalues of the community matrix defined as
$
\textbf{M}_{ij}=x_i^*\frac{\partial f_i}{\partial x_j}|_{\textbf{x}=\textbf{x}^*}
$
have negative real parts and $f_i(\textbf{x}^*)<0$ for all $i\in I_0$ \cite{goh1978sector}. Global sector stability results were also established by Goh through Lyapunov theory. Define $
U=\{\textbf{x}\text{ }|\textbf{ } x_i>0 \text{ for } i\in I_1\text{ and } x_i\geq 0 \text{ for } i\in I_0\}$. A semi-feasible equilibrium point is called globally sector stable if it is sector stable
relative to the set $U$. The semi-feasible equilibrium point $\textbf{x}^*$ of the generic population model (\ref{eq:gpm}) is globally sector stable if there  exists positive constants $d_1,d_2,\dots,d_S$ such that at every point in $U$, the function 
\begin{equation*}
    V(\textbf{x})=\sum_{i=1}^Sd_i (x_i-x_i^*)f_i(\textbf{x})\leq 0,
\end{equation*}
and it does not vanish identically along a solution of (\ref{eq:gpm}) except for $\textbf{x}=\textbf{x}^*$ \cite{goh1978sector}.
The  result establishes valuable conditions for global sector stability in the GLV model (\ref{eq: glv}): (i) there exists a positive diagonal matrix $\textbf{D}=\text{diag}(d_1,d_2,\dots,d_S)$ such that $\textbf{D}\textbf{A}+\textbf{A}^\top\textbf{D}$ is negative semi-definite (i.e., \textbf{A} is diagonally semi-stable); (ii) the expressions $r_i+\sum_{j=1}^S\textbf{A}_{ij}x_j^*\leq 0$ for all $i\in I_0$; (iii) the function $\frac{1}{2}(\textbf{x}-\textbf{x}^*)^\top (\textbf{D}\textbf{A}+\textbf{A}^\top\textbf{D})(\textbf{x}-\textbf{x}^*) + \sum_{i\in I_0}d_ix_i(r_i+\sum_{j=1}^S\textbf{A}_{ij}x_j^*)$ does not vanish identically along any solution of (\ref{eq: glv}) except for $\textbf{x}=\textbf{x}^*$. 
}

\begin{figure}[t]
    \centering
    \includegraphics[scale=0.48]{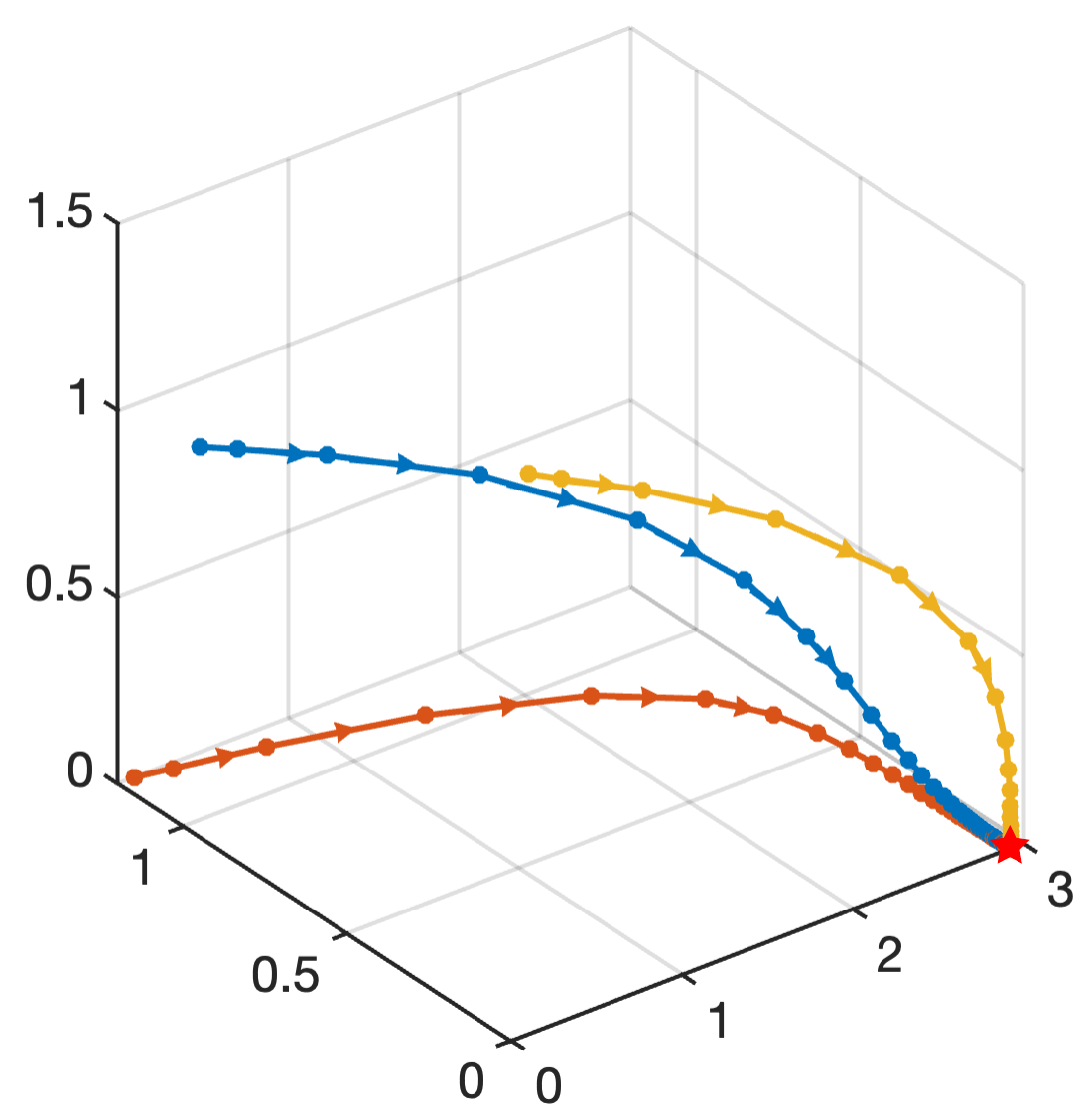}
    \caption{Trajectories of the GLV model with three initial conditions closely positioned to the other three semi-feasible equilibrium points. The globally sector stable semi-feasible equilibrium point is marked by the red star. }
    \label{fig:sector}
\end{figure}

\textcolor{black}{
However, the above global sector stability condition is difficult to verify in practice (especially when $S> 2$). Thus, Goh came up with a  conservative but simpler result. Suppose that there exists a constant matrix \textbf{G} such that 
\begin{equation}\label{eq:8.2}
\begin{split}
    &\frac{\partial f_i(\textbf{x})}{\partial x_i}|_{\textbf{x}=\textbf{x}^*}\leq \textbf{G}_{ii}<0 \quad \text{for } i=1,2,\dots,S,\\
    &\Big|\frac{\partial f_i(\textbf{x})}{\partial x_j}|_{\textbf{x}=\textbf{x}^*}\Big|\leq \textbf{G}_{ij}\quad \text{for } i\neq j
\end{split}
\end{equation}
in the set $U$. If all the leading principal minors of $-\textbf{G}$ are positive and $f_i(\textbf{x}^*)\leq 0$ for all $i\in I_0$, the semi-feasible equilibrium point $\textbf{x}^*$ is globally sector stable \cite{goh1978sector}. As an illustrative example (borrowed from \cite{goh1978sector}), consider the following GLV model:
\begin{equation*}
    \begin{cases}
        \dot{x}_1 &= x_1(11.7-4x_1-0.2x_2-0.1x_3)\\
        \dot{x}_2 &= x_2(1.2-0.8x_1-x_2-0.2x_3)\\
        \dot{x}_3 &= x_3(3-2x_1-x_2-2x_3)
    \end{cases}.
\end{equation*}
The model has four semi-feasible equilibrium points, which are $(0,1,1)$, $(0,6/5,0)$, $(0,0,3/2)$, and $(117/40,0,0)$. Let $\textbf{G}_{ii}=\textbf{A}_{ii}$ for $i=1,2,3$ and $\textbf{G}_{ij}=|\textbf{A}_{ij}|$ for $i\neq j$ such that (\ref{eq:8.2}) is satisfied. All the leading principal minors of $-\textbf{G}$ are positive. At the equilibrium point  $(117/40,0,0)$, $f_2(\textbf{x}^*)<0$ and $f_3(\textbf{x}^*)<0$. Therefore, this semi-feasible equilibrium point is globally sector stable with respect to $U=\{\textbf{x} \text{ }|\text{ } x_1>0 \text{ and } x_2,x_3\geq 0\}$, see Figure \ref{fig:sector}.
}

\textcolor{black}{
The feasibility requirement for equilibrium points in complex ecological models such as the GLV model drastically restricts their parameter space, especially when assuming random interactions. Consequently, most equilibrium points are semi-feasible, which highlights the importance and necessity of sector stability. Furthermore, Goh's results indicate that a semi-feasible equilibrium point will likely be globally sector stable if the strength of self-regulating interactions surpasses that of interspecific interactions \cite{goh1978sector}. 
}

\section{Structural stability analysis}\label{sec:9}
The notion of structural stability of dynamical systems was first introduced by Andronov and Pontryagin   under the name coarse (or rough) systems \cite{Andronov-1937}. \textcolor{black}{Different from the previous notions of stability which consider perturbations of initial conditions for a fixed dynamical system}, structural stability concerns whether the qualitative behavior of the system trajectories will be affected by small perturbations of the system model itself \cite{bavzant2000structural,peixoto1959structural,simitses2006fundamentals,thom2018structural,ziegler2013principles}. In this section, we first discuss the mathematical definition of structural stability in dynamical systems. Then, we present various metrics that can be used to quantify the structural stability of  ecological systems. 

\subsection{Mathematical definition}\label{sec:9.1}

\begin{figure}[t!]
\centering
\includegraphics[width=0.9\textwidth]{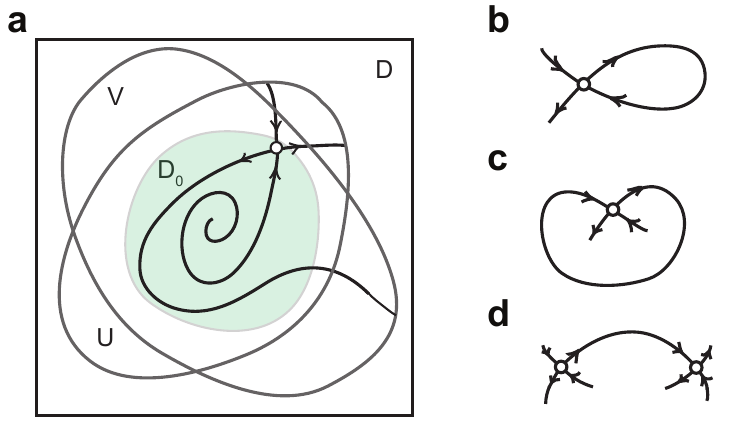}
\caption{Structural stability. (a) Andronov's definition of structural stability. (b-d) Phase portraits of structurally unstable planar systems. This figure was redrawn from \cite{Kuznetsov-Book-04}. 
}
\label{fig:structuralstability}
\end{figure}

To formally define structural stability, we introduce the concept of topologically equivalence of dynamical systems. Two dynamical systems are called topologically equivalent if there is a homeomorphism $h: \mathbb{R}^S \to \mathbb{R}^S$ mapping their phase portraits, preserving the direction of time.  Consider two smooth continuous-time dynamical systems  (i) $\dot{\textbf{x}}(t) = \textbf{f} (\textbf{x}(t))$ and (ii) $\dot{\textbf{x}}(t) = \textbf{g} (\textbf{x}(t))$. Both systems (i) and (ii) are defined in a closed region $D \in\mathbb{R}^S$, see Figure~\ref{fig:structuralstability}a. System (i) is called structurally stable in a region $D_0\subset D$ if for any system (ii) that is sufficiently $C^1$-close to system (i) there are regions $U$, $V \subset D$, and $D_0 \subset U$, $D_0 \subset V$ such that system (i) is topologically equivalent in $U$ to system (ii) in $V$, see Figure~\ref{fig:structuralstability}a.  Here, the systems (i) and (ii) are $C^1$-close if their ``distance,'' defined as
\begin{equation*}
    d = \sup_{\textbf{x}\in D} \left\{ 
    \| \textbf{f}(\textbf{x}) - \textbf{g}(\textbf{x}) \| 
+  \left\| \frac{d \textbf{f}(\textbf{x})}{d \textbf{x}} - \frac{d \textbf{g}(\textbf{x})}{d \textbf{x}} \right\|
\right\},
\end{equation*}
is small enough.

Andronov and Pontryagin  offered sufficient and necessary conditions for a two-dimensional continuous-time dynamical system to be structurally stable \cite{Andronov-1937}. A smooth dynamical system  $\dot{\textbf{x}}(t) = \textbf{f} (\textbf{x}(t))$ with $\textbf{x}(t) \in \mathbb{R}^2$,  is structurally stable in a region $D_0 \subset \mathbb{R}^2$ if and only if (i) it has a finite number of equilibrium points and limit cycles in $D_0$, and all of them are hyperbolic; (ii) there are no saddle separatrices returning to the same saddle, see Figure~\ref{fig:structuralstability}b and c, or connecting two different saddles in $D_0$, see Figure~\ref{fig:structuralstability}d.  This is often called the Andronov-Pontryagin criterion, which gives the complete description of structurally stable systems on the plane. It has been proven that a typical or generic two-dimensional system always satisfies the Andronov-Pontryagin criterion and hence is structurally stable \cite{Peixoto-1962}. In other words, structural stability is a generic property for planar systems. Yet, this is not true for high-dimensional systems. Later on, Morse and Smale established the sufficient conditions for an $S$-dimensional dynamical systems to be structurally stable \cite{Smale-1961,Smale-1967}. Such systems, often called Morse-Smale systems, have only a finite number of equilibrium points and limit cycles, all of which are hyperbolic and satisfy a transversaility condition on their stable and unstable invariant manifolds. 

\subsection{Structural stability metrics}\label{sec:9.2}
The notion of structural stability has  been  explored in diverse ecological systems. For instance, Recknagel  investigated the structural stability of aquatic ecological systems by using the catastrophe theory and applied his approach as an aid in decision-making for water quality management \cite{recknagel1985analysis}. In addition, the concept of structural stability have been heavily used in soil ecological systems \cite{garcia2012soil,perfect1990factors,perfect1990rates,sankaran1999determinants,van2007type}. 
Below, we survey various  measures that can be used to quantify the structural stability of  ecological systems.

\textcolor{black}{
Rohr et al.  introduced a mathematical framework based on the concept of structural stability to elucidate the influence of network architecture on community persistence with the GLV model \cite{rohr2014structural}. They proposed that an ecological system becomes more structurally stable as the area of the parameter space of the model expands, resulting in both a dynamically stable and feasible equilibrium. In particular, the authors investigated the range of conditions necessary for the stable coexistence of all species in mutualistic systems and showed that numerous observed mutualistic network architectures tend to maximize the volume of parameter space under which species coexist. This implies that having both a nested network architecture and a small mutualistic trade-off is one of the most preferable structures for community persistence.
}

\textcolor{black}{
Grilli et al.   developed a geometrical framework to study the range of conditions necessary for feasible coexistence of large ecological systems \cite{grilli2017feasibility}. They quantified the structural stability of an ecological system using the GLV model as the volume of the feasibility region when varying intrinsic growth rates, which can be approximated by
\begin{equation}
    \Theta \approx \Bigg(1+\frac{1}{\pi}\frac{C_A\mu_A(2\mu_d-S C\mu_A)}{\mu_d-SC^2\mu_A^2}\Bigg)^S,
\end{equation}
where $C_A$ is the connectance  of the interaction matrix \textbf{A}, $\mu_A$ is the mean of the off-diagonal elements of \textbf{A}, and $\mu_d$ is the mean of the diagonal elements of \textbf{A}. The authors further analytically predicted the range of coexistence conditions in  more than 100 empirical ecological systems. Recently, using the above approximation of structural stability, Portillo et al.  explored the correlation between structural stability and various network measures, such as centrality and modularity, and illustrated that optimal modularity has a negative impact on biological diversity (structural stability) in empirical ecological systems \cite{portillo2022global}.
}

\textcolor{black}{
During the same time, Saavedra et al.  proposed novel metrics analogous to stabilizing niche differences and fitness differences with the GLV model, which can measure the range of conditions compatible with multi-species coexistence, i.e., structural stability \cite{saavedra2017structural}.  The structural analogs of the niche and fitness differences are defined as
\begin{equation}\label{eq:feas}
        \Omega(\textbf{A})= \frac{|\text{det}(\textbf{A})|}{\sqrt[S]{\pi/2}}\int \cdots \int_{\mathbb{R}^S_{\geq 0}} e^{-\textbf{x}^\top\textbf{A}^\top\textbf{A}\textbf{x}}d\textbf{x} \quad \text{and} \quad
        \vartheta = \text{arccos}\Big(\frac{\textbf{r}^\top \textbf{r}_c}{\|\textbf{r}\|\|\textbf{r}_c\|}\Big),
\end{equation}
respectively, where $\textbf{A}$ is the interaction matrix, $\textbf{r}$ is the vector of growth rates, and $\textbf{r}_c$ is the centroid of the feasibility domain defined as 
\begin{equation*}
    \textbf{r}_c=\frac{1}{S}\Big(\frac{\textbf{a}_1}{\|\textbf{a}_1\|}+\frac{\textbf{a}_2}{\|\textbf{a}_2\|} +\cdots + \frac{\textbf{a}_S}{\|\textbf{a}_S\|}\Big).
\end{equation*}
Here, $\textbf{a}_i$ denotes the $i$th column of \textbf{A}. Therefore, feasible solutions can be fulfilled as long as \textbf{r} is inside the cone defining the domain of feasibility $D(\textbf{A})=\{\textbf{r}\in\mathbb{R}^S_{> 0} \text{ such that } \textbf{A}^{-1}\textbf{r}>0\}$. In other words, the structural analog of the fitness difference $\theta$ is small enough relative to the structural analog of the niche difference $\Omega(\textbf{A})$. The authors further applied their structural approach to a field system of annual plant competitors occurring on serpentine soils.}

\textcolor{black}{Saavedra and his coauthors have utilized similar approaches to study the structural stability of various ecological systems \cite{cenci2018structural,cenci2018rethinking,rohr2016persist,saavedra2016seasonal,saavedra2014structurally,saavedra2016nested,song2017some,song2018guideline,song2018structural}. For example, they   quantified structural stability using the quantity
\begin{equation}
    \eta = \frac{(1-\cos^2{(\vartheta)}}{\cos^2{(\vartheta)}},
\end{equation}
which measures how big the deviations are from the structural vector (geometric centroid) compatible with a positive stable equilibrium \cite{saavedra2014structurally}. They further proved that the smaller the level of global competition, the broader the conditions for having feasible solutions. In addition,  Song and Saavedra  proposed a measure of structural stability using the quantity $\Omega(\textbf{A})^{1/S}$ with a transformation $2\textbf{A}^\top\textbf{A}=\boldsymbol{\Sigma}^{-1}$, which offers an approximation to the level of external perturbations tolerated by an ecological system \cite{song2018structural}. The authors further found that their measure is the only consistent predictor of changes in species richness in empirical ecological systems among different ecological and environmental variables. 
}

\textcolor{black}{
Lately, Pettersson et al.  developed a metric called instability to quantify the proximity to collapse and level of structural stability with the GLV model, which provides deep insights into the dynamics and limits of stability and collapse of ecological systems \cite{pettersson2020predicting}. The instability metric is defined as
\begin{equation}
    \gamma(s) = \frac{\sigma_A-\sigma_f\big(S_{\text{pred}}(s)\big)}{\sigma_c\big(S_{\text{pred}}(s)\big)-\sigma_f\big(S_{\text{pred}}(s)\big)},
\end{equation}
where $\sigma_A$ is the standard derivation of interspecific interaction strengths, $\sigma_f$ is the first extinction boundary, and $\sigma_c$ is the collapse boundary. Here, $S_{\text{pred}}(s)$ is the predicted initial biodiversity in terms of the actual biodiversity $s$.
The metric $\gamma\in[0,1]$, and it can be computed from observable quantities of an ecological system. Notably, a higher $\gamma$ value indicates lower structural stability of a system, which increases the likelihood of collapse due to perturbations or external pressures.
}

\section{Higher-order stability analysis}\label{sec:10}

\begin{figure}[t]
    \centering
    \includegraphics[scale=0.4]{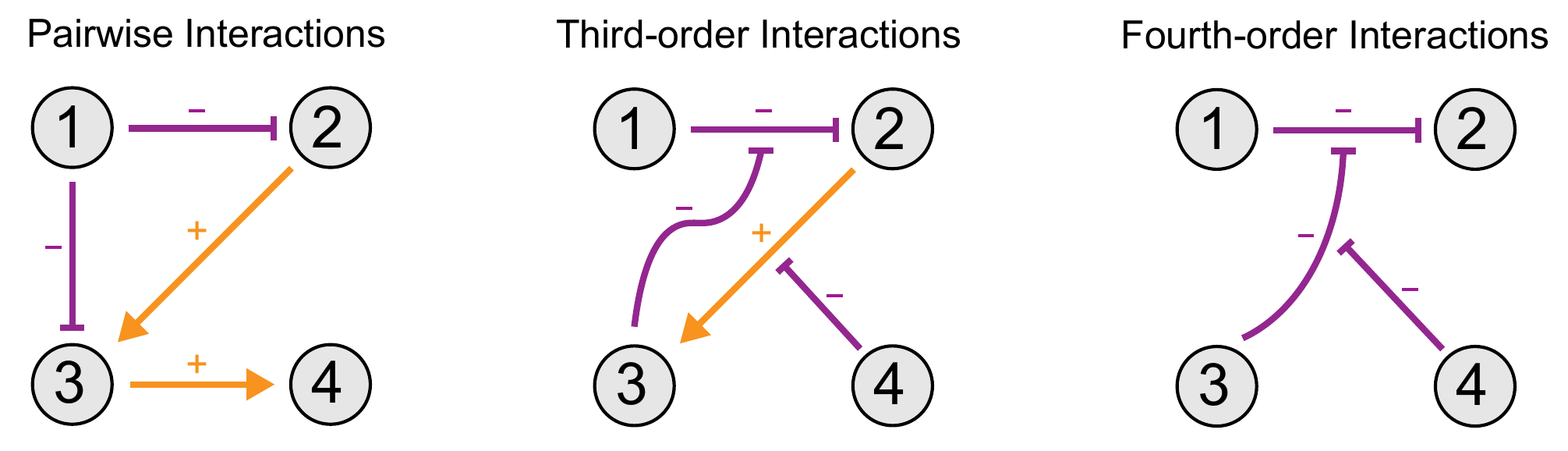}
    \caption{\textcolor{black}{Examples of pairwise, third-order, and fourth-order interactions in ecological systems. For pairwise interactions, species 1 could produce an antibiotic, inhibiting the growth of species 2. For third-order interactions,  species 3 might degrade the antibiotics produced by species 1, thus alleviating the inhibitory effect of species 1 on species 2. For fourth-order interactions, the activity of the antibiotic-degrading enzyme by species 3 may in turn be inhibited by compounds produced by species 4. This figure was redrawn from \cite{bairey2016high}.}}
    \label{fig:higher-order}
\end{figure}

Various complex networks, such as ecological networks, social networks, biological networks, and chemical reaction networks, exhibit higher-order interactions, in which the interactions go beyond pairwise relationships and involve multiple nodes or entities at the same time \cite{alvarez2021evolutionary,battiston2021physics,billick1994higher,boccaletti2023structure,chen2023survey,chen2023teasing,dotson2022deciphering,ghazanfar2020investigating,wilson1992complex}. 
In particular, in  ecological networks, species interactions frequently occur in higher-order combinations, where the relationship between two species is influenced by one or more additional species \cite{bairey2016high,kleinhesselink2022detecting,levine2017beyond}, see Figure \ref{fig:higher-order}. From the modeling perspective, higher-order interactions can be explicitly modeled by adding higher-order terms to the classical GLV model. Alternatively, they can be implicitly captured through consumer-resource models, i.e., consumer species associated with the same resource form a higher-order interaction.
While the significance of higher-order interactions has been recognized, their comprehensive impact on the stability of ecological systems has not been fully understood.

\textcolor{black}{
Pioneering research  has shown that one type of third-order interaction, in which one species mitigates the negative interactions between two others, can stabilize well-mixed ecological systems with particular network configurations \cite{kelsic2015counteraction}. In addition, increasing the order of  interactions among species can weaken the destabilizing effect and amplify the variance of species abundances at the equilibrium \cite{de2000random,yoshino2008rank}.
In this section, we first introduce the GLV model with higher-order interactions and present various results related to the stability of ecological systems with higher-order interactions. We subsequently outline two potential  approaches for analytically determining the stability of the GLV model with higher-order interactions via polynomial systems theory. Finally, we discuss various consumer-resource model for implicit higher-order interactions and their stability properties. 
} \textcolor{black}{Note that some results for higher-order stability also involve linear stability analysis.}

\subsection{GLV model with higher-order interactions}\label{sec:10.1}

\textcolor{black}{
The dynamics of ecological systems with higher-order interactions are often described by the  GLV model with higher-order interactions \cite{aladwani2019addition,singh2021higher} defined as
\begin{equation}\label{eq:hoglv}
\begin{split}
     \dot{x}_i(t) &= x_i(t)\Big[r_i + \sum_{j=1}^S \textbf{A}_{ij} x_j(t) + \sum_{j=1}^S\sum_{k=1}^S \textsf{B}_{ijk}x_j(t)x_k(t)\\ & + \sum_{j=1}^S\sum_{k=1}^S\sum_{l=1}^S \textsf{C}_{ijkl}x_j(t)x_k(t)x_l(t) + \cdots\Big]
\end{split}
\end{equation}
for $i=1,2,\dots,S$. Here, $\textsf{B}\in\mathbb{R}^{S\times S\times S}$ is a third-order tensor whose off-diagonal elements represent the effect that species $j$ and $k$ has upon species $i$.   $\textsf{C}\in\mathbb{R}^{S\times S\times S\times S}$ is a fourth-order tensor whose off-diagonal elements represent the effect that species $j$, $k$, $l$ has upon species $i$. In fact, the  GLV model with higher-order interactions (\ref{eq:hoglv}) belongs to the family of polynomial dynamical systems \cite{chen2022explicit}.  AlAdwani and Saavedra \cite{aladwani2019addition} determined the number of non-trivial equilibrium points of the GLV model with higher-order interactions (\ref{eq:hoglv}) based on Bernshtein’s theorem \cite{bernshtein1975number} from algebraic geometry. However, the stability of those non-trivial equilibrium points remains unknown. So far, most of the stability analysis of higher-order interactions depends on linearization or numerical simulations of the GLV model with higher-order interactions or other similar models.
}

Bairey et al.  found that the classical relationship between diversity and stability is reversed when considering high-order interactions  by simulating the dynamics of ecological systems using a replicator model \cite{bairey2016high}, which is defined based on the higher-order GLV model, with random interactions. Denote the bracket part in the GLV model with higher-order interactions (\ref{eq:hoglv}) as $f_i(\textbf{x})$. The replicator model is defined as
\begin{equation}
    \dot{x}_i(t) = x_i(t)\Big[f_i(\textbf{x})-\sum_{j=1}^S x_jf_j(\textbf{x})\Big],
\end{equation}
where the extra term ensures \textbf{x} to be solved in the form of relative abundances. The authors further derived a stability criterion by examining the effective pairwise interactions through linearizing the higher-order interactions around an equilibrium point: an ecological system is stable with high probability if
\begin{equation}
    \alpha + \frac{\beta}{S} + \frac{\gamma}{S^2} + \cdots < \frac{1}{S}, 
\end{equation}
where $\alpha$, $\beta$, and $\gamma$  determines the strength (variance) of the pairwise, third-order, and fourth-order interactions, respectively. Therefore, interactions of any given order destabilize systems, and the impacts of different orders are additive, scaling differently with diversity. The combination of these interactions establishes both a lower and an upper bound on  the number of species, which can be computed as
\begin{equation*}
    S_{\text{min,max}} = \frac{1-\beta\pm\sqrt{(1-\beta)^2-4\alpha\gamma}}{2\alpha}.
\end{equation*}
When the discriminant vanishes, i.e., $1-\beta=2\sqrt{\alpha\gamma}$, the two bounds narrow a defined optimal number of species $S^{*}=\sqrt{\gamma/\alpha}$.

\textcolor{black}{Later on, Grilli et al.  explored the role of higher-order interactions in competitive ecological systems \cite{grilli2017higher}. Particularly, the authors simulated different competitive ecological systems with only third-order interactions  of the form:
\begin{equation}\label{eq:third-order}
    \dot{x_i}(t) = x_i(t)\Big[\sum_{j=1}^S\sum_{k=1}^S \textsf{B}_{ijk}x_j(t)x_k(t)\Big]
\end{equation}
for $i=1,2,\dots,S$, which can be viewed as a special case of the higher-order GLV model. The interaction tensor \textsf{B} is defined from pairwise interactions and can be computed as $\textsf{B}_{ijk}=2\textbf{H}_{ij}\textbf{H}_{ik}-\textbf{H}_{ji}\textbf{H}_{jk}-\textbf{H}_{ki}\textbf{H}_{kj}$ where the first term represents the probability of species $i$ beats both species $j$ and $k$, and the remaining two terms represent the probabilities that either species $j$ or $k$ dominate. Here, the matrix \textbf{H} encodes the dominance relationships among species, i.e., $\textbf{H}_{ij}$ represents the probability that species $i$ outcompetes species $j$. Compared to the pairwise dynamics, the equilibrium point of (\ref{eq:third-order}) is unchanged, and more importantly, it is globally stable.
Therefore, they concluded that incorporating higher-order interactions into competitive ecological systems enhances stability, making species coexistence resilient to perturbations of both population abundance and parameter values \cite{grilli2017higher}.
}

\textcolor{black}{
Singh and Baruah  derived general rules for species coexistence modulated by higher-order interactions by simulating the higher-order GLV model \cite{singh2021higher}. They demonstrated that negative higher-order interactions can promote  coexistence if these interactions strengthen intraspecific competition more than interspecific competition, while positive  higher-order interactions can stabilize coexistence across a wide range of fitness differences, disregarding differences in strength of interspecific and intraspecific competition. Recently, Gibbs et al.  also exploited the higher-order GLV model and found that when interspecific higher-order interactions became excessively harmful compared to self-regulation, the stability of coexistence in diverse ecological systems will be compromised \cite{gibbs2022coexistence}. Furthermore, they found that coexistence can also be lost when these interactions are weak and mutualistic higher-order effects became dominant \cite{gibbs2022coexistence}.
}

\subsection{Two potential approaches}\label{sec:10.2}

As mentioned, the  GLV model with higher-order interactions (\ref{eq:hoglv}) belongs to the family of polynomial dynamical systems. Thus, polynomial systems theory can be leveraged to analytically study the stability properties of the model. We summarize two recent approaches that can analytically determine the stability of homogeneous polynomial dynamical systems \cite{ahmadi2019algebraic,chen2022explicit}. 
For simplicity, let's assume that a higher-order ecological system with $S$ species only contains $p$th-order interactions. The resulting  GLV model with higher-order interactions thus becomes homogeneous (of degree $p$), similar to the form (\ref{eq:third-order}). In fact, every polynomial dynamical system can be homogenized by introducing additional variables. Therefore, the two approaches can be in principle applied to analyze the stability of the GLV model with higher-order interactions.

\subsubsection{Lyapunov approach}\label{sec:10.2.1}
\textcolor{black}{
Ali and Khadir proved that the existence of a rational Lyapunov function, defined as the ratio of two polynomials, is both necessary and sufficient for asymptotic stability of a homogeneous polynomial dynamical system at the origin \cite{ahmadi2019algebraic}. The rational Lyapunov function takes the form of
\begin{equation*}
    V(\textbf{x}) = \frac{q(\textbf{x})}{(\sum_{i=1}^S x_i^2)^r},
\end{equation*}
where $r$ is a non-negative integer and $q(\textbf{x})$ is a homogeneous positive definite polynomial of degree $2r + 2$. Since the Lyapunov inequalities on both the rational function and its derivative have sum of squares certificates, $V(\textbf{x})$ can always be found by semi-definite programming \cite{vandenberghe1996semidefinite}. }

\textcolor{black}{
The authors further proved that an $S$-dimensional homogeneous polynomial dynamical system of degree $p$ is asymptotically stable (in the sense of Lyapunov) at the origin if and only if there exist a non-negative integer $r$, a positive even integer $s$, with $2r < s$, and symmetric matrices $\textbf{P} \succeq \textbf{I}$ and $\textbf{Q} \succeq \textbf{I}$, such that
\begin{equation*}
    \big\langle \textbf{z}(\textbf{x}), \textbf{Q}\textbf{z}(\textbf{x})\big\rangle = -2\|\textbf{x}\|^2 \big\langle J(\textbf{m}(\textbf{x}))^\top \textbf{P}\textbf{m}(\textbf{x}), \textbf{f}(\textbf{x})\big\rangle + 2r\textbf{m}(\textbf{x})^\top \textbf{P}\textbf{m}(\textbf{x})\big\langle \textbf{x}, \textbf{f}(\textbf{x})\big\rangle,
\end{equation*}
where $\textbf{m}(\textbf{x})$ and $\textbf{z}(\textbf{x})$ denotes the vector of monomials in \textbf{x} of degree $\frac{s}{2}$ and $\frac{s+k+1}{2}$, respectively, and $J(\textbf{m}(\textbf{x}))$ denotes the Jacobian matrix of $\textbf{m}(\textbf{x})$ \cite{ahmadi2019algebraic}. However, solving this hierarchy of semi-definite programs requires trying all possible combinations of $s$ and $r$, and the degree $s$ cannot be bounded solely based on the values of $S$ and $p$.
}

\subsubsection{Tensor decomposition approach}\label{sec:10.2.2}
\textcolor{black}{
Chen  exploited tensor algebra to determine the stability of homogeneous polynomial dynamical systems \cite{chen2022explicit}. First, the author proved that every $S$-dimensional homogeneous polynomial dynamical system of degree $p$ can be equivalently represented by a $(p+1)$th-order $S$-dimensional tensor, denoted by $\textsf{A}\in\mathbb{R}^{S\times S\times \stackrel{p+1}{\cdots}\times S}$. Suppose that \textsf{A} is orthogonally decomposable, i.e., 
\begin{equation}\label{eq:orthdecomp}
\textsf{A} = \sum_{i=1}^S \lambda_i \textbf{v}_i\circ\textbf{v}_i\circ\stackrel{p+1}{\cdots}\circ\textbf{v}_i,
\end{equation}
where $\lambda_i$ are often referred to as the Z-eigenvalues of \textsf{A} with the corresponding Z-eigenvectors $\textbf{v}_i$ \cite{robeva2016orthogonal}. 
The stability properties of the homogeneous polynomial dynamical system can be readily obtained through the Z-eigenvalues, similar to linear stability. Let the initial condition $\textbf{x}_0=\sum_{i=1}^S \alpha_i \textbf{v}_i$. The equilibrium point $\textbf{x}^{*}=\textbf{0}$ is 
    (i) stable if and only if $\lambda_i \alpha_i^{p-1}\leq 0$ for all $i$;
    (ii) asymptotically stable if and only if  $\lambda_i \alpha_i^{p-1}< 0$ for all $i$;
    (iii) unstable if and only if $\lambda_i \alpha_i^{p-1}> 0$ for some $i$.
When $p$ is odd, the initial condition can be ignored, and the stability conditions only depends on the Z-eigenvalues (exactly same as the case of  linear stability).
}

The tensor decomposition approach provides straightforward criteria for determining the stability of homogeneous polynomial dynamical systems compared to the Lyapunov approach. Notably, Chen successfully employed this approach to analytically determine the stability of an ecological system described by the GLV model with higher-order interactions (considering third-order interactions only) with supply rates  \cite{chen2022explicit}. The model is defined as
\begin{equation}
    \dot{x}_i(t) = x_i(t)\Big[\sum_{j=1}^S\sum_{k=1}^S \textsf{B}_{ijk}x_j(t)x_k(t)\Big] + s_i,
\end{equation}
where $s_i$ is the supply rate for species $i$. The stability of the model can be obtained for orthogonal decomposable \textsf{B}.
However, not all tensors can be decomposed in the form of (\ref{eq:orthdecomp}), which limits the applicability of this approach \cite{robeva2016orthogonal}. Nevertheless, we believe that both approaches (the Lyapunov and tensor decomposition approaches) have a significant potential for analytically understanding the stability of higher-order dynamics in ecological systems.

\subsection{Implicit higher-order interactions}\label{sec:10.3}
Implicit higher-order interactions have been considered in   consumer-resource models \cite{lafferty2015general}.
MacArthur's consumer-resource model is one of the most commonly used models for studying consumer-resource interactions \cite{advani2018statistical,chesson1990macarthur}. It is defined as
\begin{equation}\label{eq:cr}
    \begin{cases}
        \dot{x}_i(t) &= x_i(t)\Big[\theta_i\sum_{\alpha=1}^N\textbf{C}_{\alpha i}y_{\alpha}(t) - \mu_i\Big]\\
        \dot{y}_{\beta}(t) &= y_{\beta}(t)\Big[\rho_\beta  (1-\frac{y_\beta(t)}{K_{\beta}})-\sum_{j=1}^S \textbf{C}_{\beta j}x_j(t)\Big]
    \end{cases}
\end{equation}
for $i=1,2,\dots,S$ and $\beta=1,2,\dots, N$, where $x_i(t)$ is the abundance of species $i$, $y_\beta$ is the abundance of resource $\beta$, $\theta_i\geq 0$ is the efficiency rate of species $i$ for converting the consumed resources into biomass, $\textbf{C}_{\alpha i}$ is the rate at which species $i$ consumers resource $\alpha$, $\mu_i$ is the mortality rate for species $i$, and $\rho_\beta$ and $K_{\beta}$ are the maximal growth rate and carrying capacity of resource $\beta$, respectively.

In  MacArthur's original formulation where the consumers and resources have different timescales, it can be shown that a feasible equilibrium point implies global stability by using a minimization principle \cite{aparicio2023feasibility,arthur1969species,arthur1970species}.  Later on, the global asymptotic stability of  MacArthur's consumer-resource model with consumers and resources  evolved in the same timescale has been proved by Case and Casten \cite{case1979global}. Notably, they rewrote  MacArthur's consumer-resource model into the GLV form (\ref{eq: glv}) with 
\begin{equation*}
\footnotesize
    \tilde{\textbf{x}} = \begin{bmatrix}
        x_1\\ \vdots\\x_S\\y_1\\\vdots\\y_N
    \end{bmatrix}, \text{ } 
    \tilde{\textbf{r}}=\begin{bmatrix}
       -\mu_1\\ \vdots\\-\mu_S\\\rho_1\\\vdots\\\rho_N
    \end{bmatrix},\text{ } 
    \tilde{\textbf{A}} = \begin{bmatrix}
        0 &\cdots & 0 & \theta_1\textbf{C}_{11} & \cdots  & \theta_1\textbf{C}_{N1}\\
         \vdots & \vdots & \vdots & \vdots & \vdots & \vdots \\
        0 & \cdots & 0 & \theta_S\textbf{C}_{1S} & \cdots  & \theta_S\textbf{C}_{NS}\\
         -\textbf{C}_{11}  & \cdots & -\textbf{C}_{1S} & -\rho_1/K_1 & \cdots & 0 \\
         \vdots & \vdots & \vdots & \vdots & \vdots & \vdots \\
           -\textbf{C}_{N1}  & \cdots & -\textbf{C}_{NS} & 0 & \cdots &  -\rho_N/K_N
    \end{bmatrix},
\end{equation*}
where $\tilde{\textbf{x}}$, $\tilde{\textbf{r}}$, and $\tilde{\textbf{A}}$ are the composite vector (or matrix) of  species abundances, intrinsic growth rates, and interactions, respectively. Based on the formulated GLV model, they successfully constructed a Lyapunov function to demonstrate that any feasible equilibrium point is globally asymptotically stable.

Aparicio et al.  explored the feasibility domain of  MacArthur's consumer-resource model \cite{aparicio2023feasibility}. Based on Case and Casten's result, the feasibility domain can be automatically treated as the stability domain in  MacArthur's consumer-resource model. Specifically, the authors utilized the structural stability metrics as defined in (\ref{eq:feas}) to compute the feasibility domain. Their findings unveil that  the feasibility of  MacArthur's consumer-resource model diminishes with the pool size of the consumers increases, but the expected fraction of feasible consumers rises in this scenario \cite{aparicio2023feasibility}. Conversely, if resources exhibit linear growth, an increase in the resource pool reduces the feasibility of the model and the expected fraction of feasible consumers. However, if the resources increase logistically, the trend is reversed.

In the following, we present several results on the stability of consumer-resource models with various structure.

\subsubsection{Impacts of cascade structure}\label{sec:10.3.1}
Cascade structure plays a critical role in the stability of ecological systems \cite{shanafelt2018stability}. Schreiber investigated the global stability of consumer-resource cascades (serially arranged containers with a dynamic consumer population that receives a flow of resources from the preceding container) \cite{schreiber1996global}. The cascade model is defined as (assuming the numbers of consumers and resources are the same)
\begin{equation}\label{eq:crc}
    \begin{cases}
        \dot{x}_i(t) &= x_i(t)\Big[\theta\big(x_i(t)\big)  h(y_i(t)/x_i(t)^b)-\mu\big(x_i(t)\big)\Big]\\
        \dot{y}_i(t) & = \gamma \big(y_{i-1}(t)-y_i(t)\big)-x_i(t)h(y_i(t)/x_i(t)^b)
    \end{cases}
\end{equation}
for $i=1,2,\dots,S$, where $\theta$ is the assimilation efficiency of the consumer dependent on $x_i$, $\mu$ is the mortality rate of the consumer dependent on $x_i$, $\gamma$ is the flow rate, and $h(y_i/x_i^b)$ is the functional response for $b\in[0,1]$ \cite{arditi1991functional}. To ensure both biological realism and mathematical tractability, it is assumed that (i) $h:\mathbb{R}\rightarrow \mathbb{R}$ is $C^1$ with $h(0)=0$ and $h'(x)>0$ for all $x\geq 0$; (ii) The limit $h_{\infty}=\lim_{x\rightarrow \infty} h(x)$ exists and is finite; (iii) $\lim_{x\rightarrow \infty}xh'(x)=0$; (iv) $0<\theta\leq 1$ and $\mu>0$ for all $x_i\geq 0$; (v) $\theta'(x_i)\leq 0$ and $\mu'(x_i)\geq 0$ for all $x_i$; (vi) $\theta(0)h_{\infty}-\mu(0)>0$. Under these assumptions (i)-(vi), the author proved that there exists an equilibrium point of the model  which is globally stable. Hence, the ratio-dependent functional response can lead to persistence of the consumer population in all containers \cite{schreiber1996global}.

\subsubsection{Impacts of mutualistic interactions}\label{sec:10.3.2}

\begin{figure}[t]
    \centering
    \includegraphics[scale=0.4]{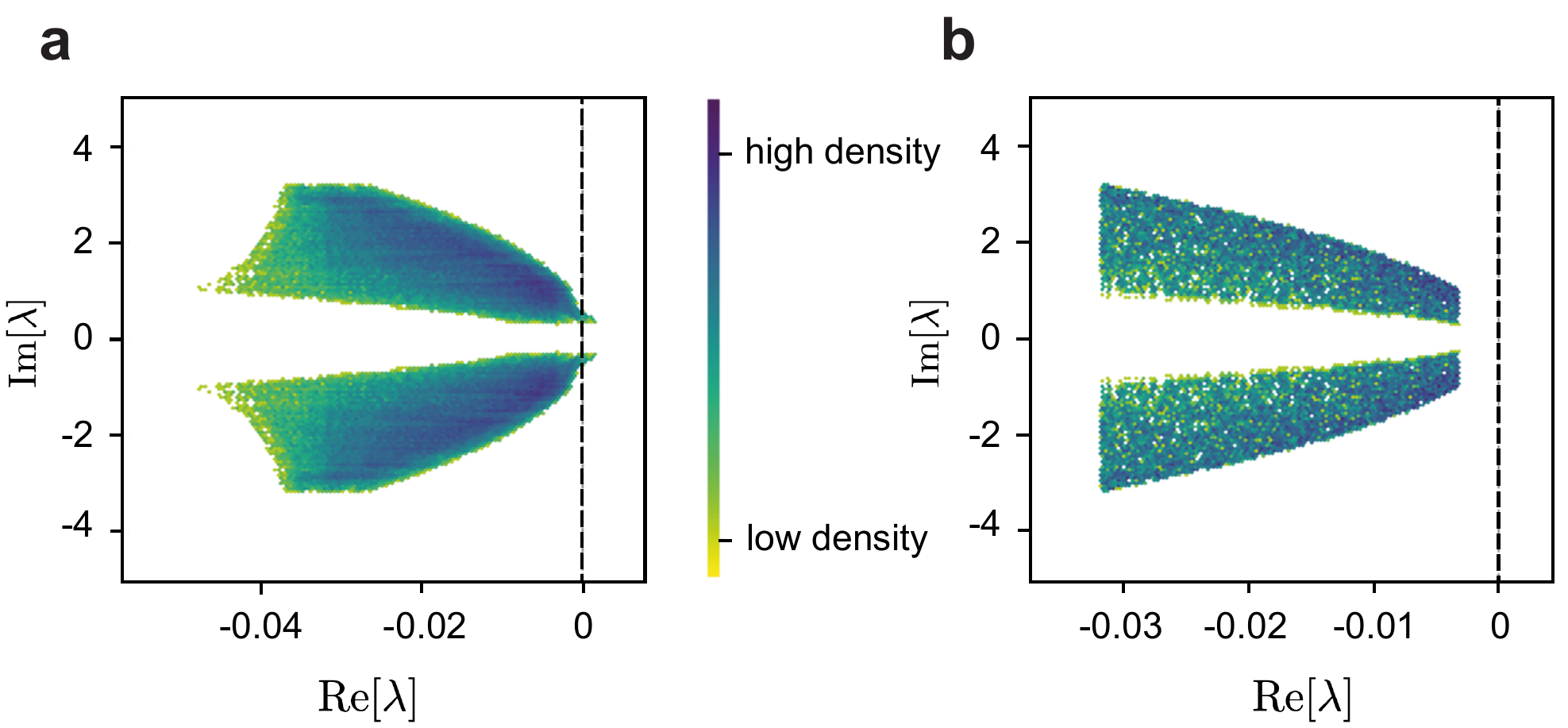}
    \caption{Eigenvalue distributions of the Jacobian matrix of the consumer-resource model with mutualistic interactions. (a) The randomly sampled production matrix \textbf{P} does not satisfy the condition (\ref{eq:crmc}), so there are positive eigenvalues. (b) The randomly sampled production matrix \textbf{P} satisfies the condition (\ref{eq:crmc}), so all eigenvalues are negative. This figure was drawn from \cite{butler2018stability}.}
    \label{fig:consumer}
\end{figure}

Microbial species capable of both resource consumption and production have the potential for mutualistic interactions \cite{coyte2015ecology,zelezniak2015metabolic}. Butler and O’Dwyer  explored the stability of a consumer-resource model with mutualistic interactions by explicitly considering the resources that microbes both consume and produce \cite{butler2018stability}. Their model is defined as (assuming the numbers of consumers and resources are the same)
\begin{equation}\label{eq:crm}
    \begin{cases}
        \dot{x}_i(t) &= x_i(t)\Big[\theta \sum_{j=1}^S \textbf{C}_{ji}y_i(t)-\sum_{j=1}^S\textbf{P}_{ji}-\mu_i\Big]\\
        \dot{y}_i(t) &= \rho_i - y_i(t)\sum_{j=1}^S\textbf{C}_{ij}x_j(t)+\sum_{j=1}^S\textbf{P}_{ij}x_j(t)
    \end{cases}
\end{equation}
for $i=1,2,\dots, S$, where $\textbf{P}_{ij}$ is the rate of the production of resource $i$ by consumer $j$, and other notations defined similarly as those in MacArthur's consumer-resource model. For simplicity, they assumed that $\textbf{C}=c\textbf{I}$ meaning each consumer specializes on a single resource. In addition, the equilibrium points for the consumers and resources are simplified as $x^*\textbf{1}$ and $y^*\textbf{1}$, where $\textbf{1}$ is the all-one vector, achieved by tuning the influx and mortality rates. This setting corresponds to the most general case of purely mutualistic interspecific interactions. The authors proved that feasible equilibrium points of the model are locally stable if  
\begin{equation}\label{eq:crmc}
    \Big(\sum_{j\neq i}\textbf{P}_{ij}\Big)^2<\frac{cx^*}{\theta}\Big(cy^*-\textbf{P}_{ii}-\frac{cx^*}{4\theta}\Big)
\end{equation}
for all $i=1,2,\dots,S$, see Figure \ref{fig:consumer}. Consequently, if microbes consume resources without producing them,  any feasible equilibrium point will be locally stable. Yet, in the presence of cross-feeding, stability is no longer assured, where it can be obtained only when mutualistic interactions are   sufficiently weak or when all pairs of taxa reciprocate each other’s assistance.

\subsubsection{Impacts of warming temperature}\label{sec:10.3.3}
Temperature can affect the stability of  ecological systems by causing changes in metabolic demands, ingestion rates, growth rates, and others \cite{amarasekare2015effects,amarasekare2019effects,donohue2016navigating,gillooly2001effects,pawar2016real}.  Synodinos et al.  used the classical Rosenzweig-MacArthur model \cite{rosenzweig1963graphical} with a Holling type II functional response to study how temperature will alter stability \cite{synodinos2023effects}. The model is defined as (assuming one consumer and one resource)
\begin{equation}
    \begin{cases}
        \dot{x}(t) &=  x(t) \Big[\frac{\theta o y(t)}{1+oly(t)}-\mu\Big]\\
        \dot{y}(t) &= y(t)\Big[\rho(1-\frac{y(t)}{K})-\frac{ox(t)}{1+oly(t)}\Big]
    \end{cases},
\end{equation}
where $o$ is the consumer attack rate, $l$ is the handling time, $\mu$ is the mortality rate of the consumer, $\rho$ is  the intrinsic growth rate of the resource, and $K$ is the carrying capacity of the resource. They analyzed the impacts of temperature on stability using the tendency for intrinsic oscillations in population abundances, which can be defined as $\varphi = olK$, i.e., the product of the three temperature-dependent parameters. The larger $\varphi$, the lower the stability of the ecological system. The authors found four different warming-stability patterns from various empirical ecological systems: stability increases, decreases, is hump-shaped or U-shaped with temperature.

\subsubsection{Impacts of dynamic switching}\label{sec:10.3.4}
The switching of the distribution of the rates of consumption between two resources is common in nature \cite{kvrivan1996optimal,kvrivan2003optimal,kvrivan1999optimal}. Gawronski et al.  introduced a dynamic switching mechanism in the classical Rosenzweig-MacArthur consumer-resource model  with one consumer and two resources \cite{gawronski2022instability}. The model is formulated as
\begin{equation}
\begin{cases}
    \dot{x}(t) &= x(t)\Big[o\frac{\theta_1 y_1(t)+\theta_2y_2(t)}{1+l\big(\theta_1 y_1(t)+\theta_2y_2(t)\big)}-\mu\Big]\\
    \dot{y}_1(t) &= y_1(t)\Big[\big(1-\textbf{L}_{11}y_1(t)-\textbf{L}_{12}y_2(t)\big)-\frac{o\theta_1x(t)}{1+l\big(\theta_1 y_1(t)+\theta_2y_2(t)\big)}\Big]\\
    \dot{y}_2(t) &= y_2(t)\Big[\big(1-\textbf{L}_{21}y_1(t)-\textbf{L}_{22}y_2(t)\big)-\frac{o\theta_2x(t)}{1+l\big(\theta_1 y_1(t)+\theta_2y_2(t)\big)}\Big]
\end{cases},
\end{equation}
where $\textbf{L}_{ij}$ represents the limitation of growth of resource $y_i$ imposed by resource $y_j$. The dynamic switching is driven by maximizing the consumer's population $x$, i.e., 
$
    \dot{\theta}_1 = (1/\tau)dx/d\theta_1,
$
where $\tau$ is the characteristic time. Through numerical simulations, the authors found that oscillations of the consumer and the mutually synchronized two resources, which emerge at $\theta_i=0.5$, become unstable for large or small $\theta_i$. Additionally, the model converges to a stable equilibrium point when either $\theta_1>0.5$ and $y_1<y_2$ or the opposite. This suggests that the consumer is unable to switch its preferred resource once it has been chosen.

\subsubsection{Impacts of periodic environments}\label{sec:10.3.5}

Periodic fluctuations and signals are widespread and strongly influence the structure and function of ecological systems \cite{bernhardt2020life,rubin2023irregular}. Bieg et al.  investigated the role of periodic environmental forcing on consumer-resource interactions \cite{bieg2023stability}. They used an extension of the classic Rosenzweig–MacArthur consumer-resource model \cite{rosenzweig1963graphical} defined as (assuming  one consumer and one resource)
\begin{equation}
    \begin{cases}
        \dot{x}(t)&= x(t)\Big[\frac{\theta o y(t)}{\alpha+y(t)}-\mu\Big]\\
        \dot{y}(t)& = y(t)\Big[\rho\big(1-\frac{y(t)}{K_{\text{mean}}+K_{\text{force}}(t)}\big)-\frac{ox(t)}{\alpha+y(t)}\Big]
    \end{cases},
\end{equation}
where $\alpha$ is the half-saturation coefficient, $K_{\text{mean}}$ is the average resource carrying capacity, and $K_{\text{force}}(t)$ is a sinusoidal function capturing periodic forces defined as
\begin{equation*}
    K_{\text{force }}(t) = A\sin{(p2\pi t)}.
\end{equation*}
Here, $A$ is  the amplitude of the  forcing, and $p$ is the forcing speed. Through local stability analysis, the authors found that the forcing speed significantly affects the stability of the model, with fast environmental forcing having a stabilizing effect. Their results further suggest that periodic fluctuations from climate change may cause sudden shifts in ecological dynamics and stability.

\subsubsection{Impacts of consumption threshold}\label{sec:10.3.6}
The idea of consumption threshold is borrowed from the basic reproduction in epidemiological models \cite{van2002reproduction}.
Duffy and Collins  introduced this notion to determine the stability of the classical Rosenzweig-MacArthur consumer-resource model \cite{duffy2016identifying}, which is defined similarly as (assuming one consumer and one resource)
\begin{equation}
    \begin{cases}
        \dot{x}(t) &= x(t)\Big[\frac{\theta o y(t)}{\alpha+y(t)}-\mu\Big]\\
        \dot{y}(t) &= y(t)\Big[\rho(1-\frac{y(t)}{K})-\frac{ox(t)}{\alpha+y(t)}\Big]
    \end{cases}.
\end{equation}
The consumption threshold is defined as $C_0=\frac{\theta oK}{\mu(\alpha+K)}$. Ecologically, it represents the parameter combination that guarantees the minimum resource consumption required for the consumer to survive. The authors proved that the trivial equilibrium point (i.e., $x^{*}=y^{*}=0$) is always unstable, and the semi-feasible equilibrium point (i.e., $x^{*}=0$ and $y^{*}=K$) is locally stable if $C_0\leq 1$. Furthermore, They also demonstrated that the model undergoes a Hopf bifurcation at
\begin{equation*}
    C_0^*=\frac{K}{\alpha+K}\Big(1+\frac{(\mu+\theta o)\alpha}{\mu K}\Big),
\end{equation*}
and the non-trivial equilibrium point (i.e., $x^*>0$ and $y^*>0$) exists for $C_0>1$ and is locally stable if $C_0<C_0^*$. In summary, when the consumption threshold $C_0\leq 1$, only the resource can survive, while when $C_0>1$, both the consumer and resource can coexist.

\subsubsection{Impacts of delayed age structure}\label{sec:10.3.7}
The introduction of delay in the consumer-resource models plays a significant role in their stability. 
Akimenko  studied the stability of a delayed age-structural consumer-resource model with $N$ resource patches and consumers  that can move freely between the patches \cite{akimenko2021stability}. The resource dynamics of the model is defined as
\begin{equation}\label{eq:dasr}
    \dot{y}_i(t) = y_i(t)\Bigg[\rho_i\Big(1-\frac{y_i(t)}{K_i}\Big) - \frac{\beta_i\hat{X}(t)}{1+\alpha_i\hat{X}(t)}\Bigg]
\end{equation}
for $i=1,2,\dots,N$, where $\beta_i>0$ is the search
rate, $\alpha_i\geq 0$ is the saturation coefficient, and $\hat{X}(t)=\int_{0}^{a_d}\kappa(a) x(a, t)da$ is the weighted quantity of
consumers at time $t$. Here, $\kappa(a)$ represents the age-specific consumer’s preference, and  $x(a,t)$ is the age-specific density of consumer population at age $a$ and time $t$. On the other hand, the consumer dynamics $x(a,t)$ is described by the delayed McKendrick-Von Foerster’s age-structured model \cite{foerster1959some,gurtin1974non,gurtin1979some} defined as
\begin{equation}\label{eq:dasc}
    \frac{\partial x}{\partial t} + \frac{\partial x}{\partial a} = -\mu\big(a, w(t-\tau)\big)x(a,t),
\end{equation}
where $\mu(a, w)$ is the age- and calorie intake rate-dependent consumer mortality rate. 
Combined with initial and boundary conditions, the author proved that the trivial equilibrium point (i.e., $y_i^*=0$ and $x^{*}(a)=0$) of the model is unstable for all $\tau>0$ and the semi-feasible equilibrium point (i.e., $y_i^*=K_i$ and $x^{*}(a)=0$) is locally asymptotically stable if the consumer’s basic reproduction number $R(w^*)<1$ and is unstable if $R(w^*)\geq 1$ for all $\tau>0$. The consumer’s basic reproduction number dependent on the calorie intake rate $w$ is defined as
\begin{equation*}
    R(w) = \int_{a_r}^{a_m} \Omega(a, w)\exp{\Big(-\int_{0}^a \mu(\zeta, w)d\zeta\Big)da},
\end{equation*}
where $\Omega(a, w)$ is the age- and calorie intake rate-dependent consumer fertility rate, and $a_r$ and $a_m$ are the age of maturation and maximum age of reproduction, respectively. Furthermore, the non-trivial equilibrium point (i.e., $y_i^*>0$ for some $i$ and $x^{*}(a)>0$) exists if and only if $R(w^*)= 1$. Define $I_0=\{i\text{ } | \text{ } y^{*}_i=0\}$. The non-trivial equilibrium point is unstable for all $\tau>0$ if $I_0\neq \emptyset$ and at least one of the following conditions is satisfied: (i) there exists $i\in I_0$ such that $\alpha_i\geq \beta_i/\rho_i$; or (ii) $\alpha_i<\beta_i/\rho_i$ for all $i\in I_0$ and $\hat{X}^*\leq \max_{i\in I_0}(\beta_i/\rho_i-\alpha_i)^{-1}$. It is locally asymptotically stable for all $\tau>0$ if one of the following conditions is satisfied: (iii) $I_0\neq \emptyset$, $\alpha_i<\beta_i/\rho_i$ for all $i\in I_0$, and $\hat{X}^*>\max_{i\in I_0}(\beta_i/\rho_i-\alpha_i)^{-1}$; (iv) $I_0= \emptyset$. Akimenko's findings  show that the digestion period (i.e., the delay parameter)  $\tau$ does not destabilize the consumer population at the trivial, semi-feasible or non-trivial equilibrium points.

Recently, El-Doma  improved the stability criteria proposed by Akimenko in analyzing the delayed age-structural consumer-resource model (\ref{eq:dasr}) and (\ref{eq:dasc}) \cite{el2023age}. Notably, the author proved that the non-trivial equilibrium point (with $I_0\neq \emptyset$) of the model is locally asymptotically stable if $X^*(\beta_i/\rho_i-\alpha_i) > 1$ and is unstable if $X^*(\beta_i/\rho_i-\alpha_i) < 1$ for $i\in I_0$. Clearly, the conditions are simpler than Akimenko's.

\subsubsection{A case study of larch moth interactions}\label{sec:10.3.8}
Plant's quality plays an important role in the population dynamics of consumers \cite{jang2009dynamics,jang2009models}. 
Din and Khan explored the dynamics between budmoths and the quality of larch trees located in the Swiss Alps \cite{din2021discrete}. They proposed an improved discrete-time consumer-resource model to capture the interactions between the plant quality index $x$ and the moth population density $y$ defined as
\begin{equation}
    \begin{cases}
        x(t+1) & = (1-v)+vx(t)-\frac{u y(t)}{\alpha+y(t)}\\
        y(t+1) &= y(t)\exp{\Big(\rho(\frac{x(t)}{\beta + x(t)}-\frac{y(t)}{K})\Big)}
    \end{cases},
\end{equation}
where $0<v<1$ is the plant vulnerability, $u$ is the moth's maximal uptake rate of the plant, $\alpha$ and $\beta$ are the half-saturation coefficients, $\rho$ is the growth rate of  moths, and $K$ is the carrying capacity of moths. The authors proved that the model has a non-trivial equilibrium point, and it is locally asymptotically stable if and only if
\begin{equation}
\begin{split}
    \frac{\beta\alpha\rho u y^*}{(\beta+x^*)^2(\alpha+y^*)^2}+\frac{(v+1)(2K-\rho y^*)}{K}&>0,\\
    v+\frac{\beta\alpha\rho u y^*}{(\beta+x^*)^2(\alpha+y^*)^2} - \frac{v\rho y^*}{K}&<1.
\end{split}
\end{equation}
The model also undergoes a period-doubling bifurcation at the non-trivial equilibrium point. Significantly, Din and Khan validated their theoretical findings with experimental data.

\section{Discussion}\label{sec:11}
\textcolor{black}{
The seven stability notions we reviewed, including linear stability, sign stability, diagonal stability, D-stability, total stability, structural stability, and higher-order stability, are not mutually exclusive but rather complement each other,  which offers a multifaceted approach to explore the stability of  ecological systems. Next, we discuss the pros and cons of each stability notion and outline its potential future prospects. 
}

\textcolor{black}{
Linear stability analysis is the most widely used approach to study the stability of  ecological systems due to its simplicity. It has been heavily used in ecological systems with various network structures, as presented in Section \ref{sec:3}. However, it cannot precisely capture the global dynamics and complex behaviors of  ecological systems that are highly nonlinear in reality. Nevertheless, it will be interesting to perform the linear stability analysis on more complex models, such as the GLV model with higher-order interactions. For example, how to express the community matrix \textbf{M} in terms of the interaction matrix \textbf{A} and the higher-order interaction tensors, and how to establish stability conditions based on higher-order interaction strengths and correlations with various higher-order network structures?
}

\textcolor{black}{
Sign stability analysis focuses on the qualitative nature of interactions in  ecological systems, disregarding of their magnitudes. It is particularly advantageous when quantitative data on interaction strengths are missing. However, the conditions described in  Section \ref{sec:4} are too stringent to realize in real ecological systems. Hence, it will be intriguing to explore simpler and more realizable conditions that are associated to the underlying network structure for the sign stability of an ecological system. Generalizing the notion to  nonlinear models, e.g., the GLV model or the GLV model with higher-order interactions, also presents exciting avenues for further work.
}

\textcolor{black}{
Diagonal stability, D-stability, and total stability analyses contribute to enhancing  the understanding of the stability of the GLV model  through characterizing the interaction matrix \textbf{A}. Moreover, they provide an insightful exploration of how the network structure encoded in the interaction matrix impact the overall stability of an ecological system. Nevertheless, characterizing diagonally stable, D-stable, or totally stable matrices is increasingly difficult for high-dimensional systems. Therefore, it will be worthwhile to further explore various network structures that can be used to establish necessary or sufficient conditions for achieving diagonally stability, D-stability, or total stability with the GLV model. Extending the current results to the GLV model with higher-order interactions also holds promising avenues for future research.}

\textcolor{black}{
Sector stability analysis is particularly useful when dealing with semi-feasible equilibrium points, which are prevalent in modeling ecological systems. Local sector stability can be readily assessed, but a definitive algorithmic approach for verifying global sector stability remains elusive. For example, how can we efficiently compute the crucial constant matrix \textbf{G} defined in Section \ref{sec:8}? Moreover, it will be worthwhile to explore simplified conditions on the sector stability of other nonlinear ecological  models such as the GLV model with higher-order interactions and MacArthur’s consumer-resource model.
}

\textcolor{black}{
Structural stability analysis assesses how the overall behavior of an ecological system is maintained under small perturbations, which is different from all the previous stability notions. Structural stability can offer unique insights into the feasibility and robustness of  ecological systems. An interesting avenue for future research would involve comparing various structural stability metrics and understanding their implications for the underlying  network structure of an ecological system.  Nonetheless, all the current structural stability metrics are developed based on the GLV model. Thus, extending these metrics to more complicated models, such as the GLV model with higher-order interactions, would also be a valuable exploration. 
}

\textcolor{black}{
Higher-order stability analysis recently becomes popular due to the importance of higher-order interactions in  ecological systems. It takes into account the nonlinearity inherent in ecological systems, providing more accurate representations, e.g., the GLV model with higher-order interactions. Yet,  the stability analysis of such higher-order models is challenging. As interactions become more complex, involving multiple species simultaneously, the mathematical complexity of the analysis increases significantly. As mentioned in  Section \ref{sec:10}, the Lyapunov and tensor decomposition approaches are promising for investigating the stability of the GLV model with higher-order interactions, which often possesses multiple feasible equilibrium points. 
}

\section{Conclusion}\label{sec:12}
\textcolor{black}{
In this article, we  presented a comprehensive and systematic review of the literature concerning the stability of  ecological systems. Notably, we inclusively surveyed various stability notions that arise in ecological systems, encompassing linear stability, sign stability, diagonal stability, D-stability, total stability, sector stability, structural stability, and higher-order stability.}

\textcolor{black}{
As mentioned earlier, empirically applying this theory to real ecological systems poses a significant challenge. First, most of the stability criteria are derived from the random matrix theory, which cannot accurately capture real ecological systems. For instance, the food web topology studies have demonstrated the non-randomness of ecological  systems, rendering to incorporate more realistic network structures \cite{landi2018complexity}. Second, the dynamics governing real ecological systems could not adhere to the GLV model, or even unknown in many cases. Thus, the stability analysis of real ecological systems may not align with the predictions of the theory. 
}

\begin{table}[t]
\caption{Summary of the most commonly used databases for real ecological systems. Links to the databases are provided.}
\label{tab:data}
\begin{tabularx}{\textwidth}{|T|X|T|}
  \hline                   
Database  &  Description & References  \\
\hline
\href{https://www.web-of-life.es}{Web of Life}                   & A graphical user interface for visualizing and downloading data on plant-animal mutualistic  ecological networks.  & \cite{fortuna2014web}  \\ \hline
\href{https://networkrepository.com/eco.php}{Network \newline Repository}                  & A network repository containing hundreds of real-world networks including ecological networks.   & \cite{rossi2015the}  \\ \hline
\href{https://www.globalbioticinteractions.org}{GloBI}                  & An extensible and open-source infrastructure for importing, searching, and exporting species-interaction data.  & \cite{poelen2014global}  \\ \hline
\href{https://www.neonscience.org}{NEON}                & Monitoring ecological systems, including freshwater and terrestrial systems across the United States. & \cite{https://doi.org/10.48443/tx5f-dy17,https://doi.org/10.48443/zayg-n802}   \\ \hline
\href{https://canberra.edu.au/globalwebdb/}{Global Web}                & An online collection of food webs. & N/A  \\ \hline
\href{https://www.sussex.ac.uk/lifesci/ebe/dopi/about\#:~:text=The\%20Database\%20of\%20Pollinator\%20Interactions\%20(DoPI)\%20documents\%20British\%20pollinator\%2D,into\%20a\%20single\%20online\%20depository}{DoPI}                & Comprising ecological networks of British pollinator-plant interactions from the published scientific literature or submitted datasets.
&  \cite{balfour2018british,goulson2015bee,ollerton2014extinctions} \\ \hline
\end{tabularx}
\end{table}

\textcolor{black}{
Additionally, applying those various stability analysis to empirical ecological systems requires accurate estimation of the species interaction matrix or community matrix. For example, the trophic flows between consumers and resources in a food web with eco-path model can be translated into the interaction coefficients of the GLV model \cite{jacquet2016no}. However, inferring the interaction network for large ecological systems, e.g., microbial communities, is challenging. Alternative frameworks, such as sensitivity testing of species abundance in response to the presence of additional species \cite{yonatan2022complexity}, may prove valuable for assessing the connectance and interaction strength of the interaction matrix or community matrix. To facilitate research in this field, we have compiled a summary of the most commonly used databases for real ecological systems in Table \ref{tab:data}.
}

\textcolor{black}{
Finally, these stability analyses could be applied to other systems beyond ecological systems, such as power grid and gene regulatory networks.  For example, the stability of power-grid networks in the homogeneous state, where the voltage frequencies of generators are all equal to a constant reference frequency, requires  the largest eigenvalues of the Jacobian matrix to be negative \cite{molnar2021asymmetry}.  Similarly, the stability of gene regulatory networks in regimes where the lifetime of mRNAs is much shorter than that of proteins depends solely on the eigenvalues of the Jacobian matrix. Furthermore, in the case of random gene regulatory networks, the underlying system becomes unstable when either the system size exceeds a critical threshold or the regulation strength becomes too strong, as determined using  random matrix theory, which aligns with  May's theory \cite{chen2002stability}.
}

\section*{Acknowledgements}

This work was supported by the National Institutes of Health R01AI141529.

\section*{Declaration of Generative AI and AI-assisted technologies in the writing process}
During the preparation of this work, the authors used ChatGPT in order to improve readability and language. After using this tool, the authors reviewed and edited the content as needed and take full responsibility for the content of the publication.

\bibliographystyle{plain}
\bibliography{references.bib}
\end{document}